\documentclass[11pt,twoside,leqno]{aomamlt2e} 
  \pageno{819}
\received{June 30, 2002}

\newtheorem{theorem}{Theorem}[section]
\newtheorem{Lemma}[theorem]{Lemma}
\newtheorem{prop}[theorem]{Proposition}
\newtheorem{Cor}[theorem]{Corollary}
\theorembodyfont{\upshape}
\newtheorem{defi}{{\it Definition}}[section]
\newtheorem{Example}{{\it Example}}[section]
\newtheorem{remark}{{\it Remark}}[section]
 
\newcommand\pref[1]{(\ref{#1})}
\renewcommand{\theequation}{\thesection.\arabic{equation}} 

\def\ca #1{{\cal #1}}
\font\gto=eufm10
\def\gt #1{\hbox{\gto #1}}
 \begin{document}
\currannalsline{159}{2004}

\let\de=\partial
\let\eps=\varepsilon
\let\phe=\varphi
\def\Hol{\mathop{\rm Hol}\nolimits}
\def\Hom{\mathop{\rm Hom}\nolimits}
\def\Fix{\mathop{\rm Fix}\nolimits}
\def\Sing{\mathop{\rm Sing}\nolimits}
\def\End{\mathop{\hbox{\rm End}}}
\def\Ker{\mathop{\rm Ker}\nolimits}
\def\Res{\mathop{\hbox{\rm Res}}}
\def\Re{\mathop{\rm Re}\nolimits}
\def\Im{\mathop{\rm Im}\nolimits}
\def\sp{\mathop{\rm sp}\nolimits}
\def\id{\mathop{\rm id}\nolimits}
\def\rk{\mathop{\rm rk}\nolimits}
\def\hz{\hat{z}}
\def\hf{\hat{f}}
\def\hg{\hat{g}}
\def\ha{\hat{a}}
\def\hb{\hat{b}}
\mathchardef\void="083F
\mathchardef\ellb="0960
\mathchardef\taub="091C
\def\C{\mathbb{C}}
\def\N{\mathbb{N}}
\def\P{\mathbb{P}}
\def\Q{\mathbb{Q}}
\def\Downarrow{\Big\downarrow}
\def\mapright#1{\smash{\mathop{\longrightarrow}\limits^{#1}}}
\def\rmapdown#1{\Big\downarrow\rlap{$\vcenter{\hbox{$\scriptstyle#1$}}$}}
\def\lmapdown#1{\llap{$\vcenter{\hbox{$\scriptstyle#1$}}$}\Big\downarrow}
\def\arrowbar{\top}
\def\Downmapsto{\vcenter{\offinterlineskip\halign{\hfil$##$\hfil\cr
\arrowbar\cr\Downarrow\cr}}}

 \title{Index theorems for holomorphic self-maps}

 \acknowledgements{}
\twoauthors{Marco Abate, Filippo Bracci,}{Francesca Tovena} 

 \institution{ Universit\`a di Pisa, Pisa, Italy\\
\email{abate@dm.unipi.it}\\
\vglue-9pt
 Universit\`a di Roma
``Tor Vergata'', Roma, Italy\\
{\scriptsize{\it E-mail addresses}\/: {\rm fbracci@mat.uniroma2.it}}\\
\hglue70pt
{\rm tovena@mat.uniroma2.it}}


  \vfil
\centerline{\bf Introduction}
\vglue8pt

The usual index theorems for holomorphic self-maps, like for instance the
classical
holomorphic
Lefschetz theorem (see, e.g., [GH]), assume that the fixed-points set
contains only
isolated points.
The aim of this paper, on the contrary, is to prove index theorems
for holomorphic self-maps having a positive dimensional fixed-points set.

The origin of our interest in this problem lies in holomorphic dynamics. A
main tool for
the complete generalization to two complex variables of the classical
Leau-Fatou flower
theorem
for maps tangent to the identity achieved in [A2] was an index theorem for
holomorphic
self-maps of
a complex surface fixing pointwise a smooth complex curve~$S$. This theorem
(later
generalized
in~[BT] to the case of a singular~$S$) presented uncanny similarities with
the Camacho-Sad
index
theorem for invariant leaves of a holomorphic foliation on a complex surface
(see~[CS]). So
we
started to investigate the reasons for these similarities; and this paper
contains what we
have
found.

The main idea is that the simple fact of being pointwise fixed by
a holomorphic self-map~$f$ induces a lot of structure on a
(possibly singular) subvariety~$S$ of a complex manifold~$M$.
First of all, we shall introduce (in \S 3) a canonically
defined holomorphic section~$X_f$ of the
bundle~$TM|_S\otimes(N_S^*)^{\otimes\nu_f}$, where $N_S$ is the
normal bundle of~$S$ in~$M$ (here we are assuming $S$ smooth;
however, we can also define $X_f$ as a section of a suitable sheaf
even  when $S$ is singular --- see Remark 3.3 --- but it turns out
that only the behavior on the regular part of $S$ is relevant for
our index theorems), and $\nu_f$ is a positive integer, the {\it order of contact} of~$f$ with~$S$, measuring how close $f$ is to
being the identity in a neighborhood~$S$ (see \S 1). Roughly
speaking, the section $X_f$ describes the directions in which $S$
is pushed by~$f$; see Proposition~8.1 for a  more precise
description of this phenomenon when $S$ is a hypersurface.

The canonical section $X_f$ can also be seen as a morphism from
$N_S^{\otimes\nu_f}$ to~$TM|_S$; its image $\Xi_f$ is the {\it canonical distribution.} When $\Xi_f$ is contained in~$TS$ (we
shall say that $f$ is {\it tangential\/}) and integrable (this
happens for instance if $S$ is a hypersurface), then on $S$ we get
a singular holomorphic foliation induced by $f$ --- and this is a
first concrete connection between our discrete dynamical theory
and the continuous dynamics studied in foliation theory. We
stress, however, that we get a well-defined foliation on $S$ {\it
only,} while in the continuous setting one usually assumes that
$S$ is invariant under a foliation defined in a {\it whole
neighborhood} of~$S$. Thus even in the tangential codimension-one
case our results will not be a direct consequence of foliation
theory.

As we shall momentarily discuss, to get index theorems it is important to
have a section of
$TS\otimes(N_S^*)^{\otimes\nu_f}$ (as in the case when $f$ is tangential)
instead of
merely a
section of $TM|_S\otimes(N_S^*)^{\otimes\nu_f}$. In Section~3, when $f$ is
not tangential
(which is a
situation akin to dicriticality for foliations; see Propositions~1.4
and~8.1)  we shall
define
other holomorphic sections~$H_{\sigma,f}$ and~$H_{\sigma,f}^1$ of
$TS\otimes(N_S^*)^{\otimes\nu_f}$
which are as good as~$X_f$ when~$S$ satisfies a geometric condition which we
call  {\it
comfortably embedded} in~$M$, meaning, roughly speaking, that $S$ is a
first-order
approximation of
the zero section of a vector bundle (see \S 2 for the precise
definition, amounting to
the
vanishing of two sheaf cohomology classes --- or, in other terms, to the
triviality of two
canonical
extensions of~$N_S$).

The canonical section is not the only object we are able to associate
to~$S$. Having a
section~$X$
of~$TS\otimes F^*$, where $F$ is any vector bundle on~$S$, is equivalent to
having an
$F^*$-valued
derivation~$X^\#$ of the sheaf of holomorphic functions~$\ca O_S$ (see
\S 5). If $E$
is another
vector bundle on~$S$, a {\it holomorphic action of~$F$ on~$E$ along~$X$}
is a ${\Bbb C}$-linear map
$\tilde
X\colon\ca E\to\ca F^*\otimes\ca E$ (where $\ca E$ and $\ca F$ are the
sheafs of germs of
holomorphic
sections of~$E$ and~$F$) satisfying $\tilde X(gs)=X^\#(g)\otimes s+g\tilde
X(s)$ for any
$g\in\ca
O_S$ and $s\in E$; this is a generalization of the notion of
$(1,0)$-connection on~$E$ (see
Example~5.1).

In Section~5 we shall show that when $S$ is a hypersurface and $f$ is
tangential (or $S$ is
comfortably embedded in~$M$) there is a natural way to define a holomorphic
action of
$N_S^{\otimes\nu_f}$ on~$N_S$ along~$X_f$ (or along~$H_{\sigma,f}$
or~$H_{\sigma,f}^1$).
And this
will allow us to bring into play the general theory developed by Lehmann and
Suwa (see,
e.g., [Su])
on a cohomological approach to index theorems. So, exactly as Lehmann and
Suwa generalized,
to any
dimension, the Camacho-Sad index theorem, we are able to generalize the index
theorems of
[A2] and
[BT] in the following form (see Theorem~6.2):

\begin{theorem}
Let $S$ be a compact{\rm ,} globally irreducible{\rm ,} possibly
singular
hypersurface in an $n$\/{\rm -}\/dimensional complex manifold $M$. Let $f\colon M
\to M${\rm ,}
$f\not\equiv \id_M${\rm ,} be a holomorphic self-map of $M$ fixing pointwise
$S${\rm ,} and denote by
$\Sing(f)$ the zero set of~$X_f$.
Assume that
\begin{itemize}
\ritem{\rm (a)} $f$ is tangential to $S${\rm ,} and then set $X=X_f${\rm ,} or that
\ritem{\rm(b)} $S^0=S\setminus\bigl(\Sing(S)\cup
\Sing(f)\bigr)$ is comfortably embedded into $M${\rm ,} and then set
$X=H_{\sigma,f}$ if
$\nu_f>1${\rm ,} or
$X=H_{\sigma,f}^1$ if $\nu_f=1$.
\end{itemize}
\noindent Assume moreover $X \not\equiv O$ {\rm (}\/a condition always satisfied when
$f$ is
tangential\/{\rm ),} and
denote by $\Sing(X)$ the zero set of~$X$. Let
$\Sing(S)\cup\Sing(X)=\bigcup_\lambda \Sigma_\lambda$ be the decomposition
of
$\Sing(S)\cup\Sing(X)$ in connected components. Finally{\rm ,} let $[S]$ be the
line bundle on~$M$
associated to the divisor~$S$. Then there exist complex numbers~$\hbox{\rm
Res}(X,S,
\Sigma_\lambda)\in\C$ depending only on the local behavior of~$X$ and $[S]$
near~$\Sigma_\lambda$
such that
$$
\sum_{\lambda}\hbox{\rm Res}(X, S,\Sigma_\lambda)=\int_S c_1^{n-1}([S]),
$$
where $c_1([S])$ is the first Chern class of $[S]$.
\end{theorem}

Furthermore, when $\Sigma_\lambda$ is an isolated point~$\{x_\lambda\}$, we
have explicit
formulas
for the computation of the residues~$\hbox{\rm Res}(X, S,\{x_\lambda\})$;
see Theorem~6.5. 

Since $X$ is a global section of
$TS\otimes(N_S^*)^{\otimes\nu_f}$, if $S$ is smooth and $X$ has
only isolated zeroes it is well-known that the top Chern
class $c_{n-1}\bigl(TS\otimes(N_S^*)^{\otimes\nu_f}\bigr)$ counts
the zeroes of~$X$. Our result shows that $c_1^{n-1}(N_S)$
is related in a similar (but deeper) way to the zero set
of~$X$. See also Section~8 for examples of results one can obtain
using both Chern classes together.

If the codimension of $S$ is greater than one, and $S$ is smooth, we can
blow-up $M$
along~$S$; then
the exceptional divisor~$E_S$ is a hypersurface, and we can apply to it the
previous
theorem. In
this way we get (see Theorem~7.2):

\begin{theorem} Let $S$ be a compact complex submanifold of
codimension $1<m<n$ in
an $n$\/{\rm -}\/dimensional complex manifold $M$. Let $f\colon M
\to M${\rm ,} $f\not\equiv \id_M${\rm ,} be a holomorphic self\/{\rm -}\/map of $M$ fixing pointwise
$S${\rm ,} and assume that $f$ is tangential{\rm ,} or that
$\nu_f>1$ {\rm (}\/or both\/{\rm ).} Let $\bigcup_\lambda \Sigma_\lambda$ be the
decomposition in connected
components of the  set of
singular directions {\rm (}\/see {\rm \S  7} for the definition\/{\rm )} for~$f$ in~$E_S$.
Then there exist
complex
numbers~$\hbox{\rm Res}(f,S,
\Sigma_\lambda)\in\C${\rm ,} depending only on the local
behavior of~$f$ and~$S$ near~$\Sigma_\lambda${\rm ,} such that
$$
\sum_{\lambda}\hbox{\rm Res}(f, S,\Sigma_\lambda)=\int_S\pi_*
c_1^{n-1}([E_S]),
$$
where $\pi_*$ denotes   integration along the fibers of the bundle $E_S\to
S$.
\end{theorem}

Theorems 0.1 and 0.2  are only two of the index theorems we can
derive using this
approach. Indeed, we are also able to obtain versions for holomorphic
self-maps
of two other main index theorems of foliation theory, the Baum-Bott index
theorem and the
Lehmann-Suwa-Khanedani (or variation) index theorem: see Theorems~6.3, 6.4,
6.6, 7.3
and~7.4.
In other words, it turns out that the existence of holomorphic actions of
suitable complex vector bundles defined only on $S$ is an efficient tool to
get
index  theorems, both in our setting  and (under slightly different
assumptions) in foliation theory; and this is another reason for the
similarities
noticed in~[A2].

Finally, in Section~8 we shall present a couple of applications of our
results to the
discrete dynamics of holomorphic self-maps of complex surfaces, thus closing
the circle and
coming
back to the arguments that originally inspired our work.

\section{The order of contact}

Let $M$ be an $n$-dimensional complex manifold, and $S\subset M$
an irreducible subvariety of codimension~$m$. We shall denote by
$\ca O_M$ the sheaf of holomorphic functions on~$M$, and by $\ca
I_S$ the subsheaf of functions vanishing on~$S$. With a slight
abuse of notations, we shall use the same symbol to denote both a
germ at $p$ and any representative defined in a neighborhood
of~$p$.
We shall denote by $TM$ the holomorphic tangent bundle of~$M$, and by $\ca
T_M$ the sheaf of
germs of local holomorphic sections of~$TM$.
Finally, we shall denote by $\End(M,S)$ the set of (germs about $S$ of)
holomorphic
self-maps of~$M$ fixing~$S$ pointwise.

Let $f\in\End(M,S)$ be given, $f\not\equiv\id_M$,
and take $p\in S$. For every
$h\in\ca O_{M,p}$ the germ $h\circ f$ is well-defined, and we have $h\circ
f-h\in\ca I_{S,p}$.

\begin{defi} The {\it $f$-order of vanishing at~$p$ of} $h\in\ca
O_{M,p}$ is given by
$$
\nu_f(h;p)=\max\{\mu\in\N\mid h\circ f-h\in\ca I_{S,p}^\mu\},
$$
and the {\it order of contact $\nu_f(p)$ of} $f$ {\it at} $p$ with~$S$ by
$$
\nu_f(p)=\min\{\nu_f(h;p)\mid h\in\ca O_{M,p}\}.
$$
\end{defi}

We shall momentarily prove that $\nu_f(p)$ does not depend on $p$.

Let $(z^1,\ldots,z^n)$ be local coordinates in a neighborhood
of~$p$. If $h$ is any holomorphic function defined in a
neighborhood of~$p$,  the definition of
order of contact yields the important relation
\begin{equation}
h\circ f- h=\sum_{j=1}^n(f^j-z^j)\,{\de h\over\de z^j}\qquad
\hbox{(mod $\ca I_{S,p}^{2\nu_f(p)}$)},
\label{eqbasedue}
\end{equation}
where $f^j=z^j \circ f$. 

As a consequence we have

\begin{Lemma} \label{fattouno} {\rm (i)} Let $(z^1,\ldots,z^n)$ be any set of local
coordinates
at~$p\in S$. Then
$$
\nu_f(p)=\min_{j=1,\ldots,n}\{\nu_f(z^j;p)\}.
$$
\begin{itemize}
\item[{\rm (ii)}] For any $h\in\ca O_{M,p}$ the function $p\mapsto\nu_f(h;p)$
is
constant in a neighborhood of~$p$.
\item[{\rm (iii)}] The function $p\mapsto\nu_f(p)$ is constant.
\end{itemize}
\end{Lemma}

\Proof  (i) Clearly, $\nu_f(p)\le\min_{j=1,\ldots,n}\{\nu_f(z^j;p)\}$. The
opposite
inequality follows from \pref{eqbasedue}.
\vglue4pt
(ii) Let $h\in\ca O_{M,p}$, and choose a set $\{\ell^1,\ldots,\ell^k\}$ of
generators
of~$\ca I_{S,p}$. Then we can write
\begin{equation}
h\circ f- h=\sum_{|I|=\nu_f(h;p)}\ell^Ig_I,
\label{eqlocal}
\end{equation}
where $I=(i_1,\ldots,i_k)\in\N^k$ is a $k$-multi-index,
$|I|=i_1+\cdots+i_k$, $\ell^I=(\ell^1)^{i_1}\cdots(\ell^k)^{i_k}$
and $g_I\in\ca O_{M,p}$. Furthermore, there is a multi-index $I_0$
such that $g_{I_0}\notin\ca I_{S,p}$. By the coherence of the
sheaf of ideals of $S$, the relation \pref{eqlocal}  holds for the
corresponding germs at all points~$q\in S$ in a neighborhood
of~$p$. Furthermore, $g_{I_0}\notin\ca I_{S,p}$ means that
$g_{I_0}|_S\not\equiv 0$ in a neighborhood of $p$, and thus
$g_{I_0}\notin\ca I_{S,q}$ for all $q\in S$ close enough to~$p$.
Putting these two observations together we get the assertion.
\vglue4pt
(iii) By (i) and (ii) we see that the function
$p\mapsto\nu_f(p)$ is locally constant and since  $S$ is connected, it is
constant everywhere.\Endproof 
\vglue4pt

We shall then denote by $\nu_f$ the {\it order of contact} of $f$ with~$S$,
computed
at any point~$p\in S$. 

As we shall see, it is important to compare the order of contact of~$f$ with
the
$f$-order of vanishing of germs in~$\ca I_{S,p}$.

\begin{defi}  We say that $f$ is {\it
tangential} at~$p$ if
$$
\min\bigl\{\nu_f(h;p)\mid h\in\ca I_{S,p}\bigr\}>\nu_f.
$$
\end{defi}

\begin{Lemma}  \label{fattodue} Let $\{\ell^1,\ldots,\ell^k\}$ be a set of
generators of~$\ca
I_{S,p}$. Then
$$
\nu_f(h;p)\ge\min\{\nu_f(\ell^1;p),\ldots,\nu_f(\ell^k;p),\nu_f+1\} 
$$
for all  $h\in\ca I_{S,p}$.
In particular{\rm ,} $f$ is tangential at~$p$ if and only if
$$
\min\{\nu_f(\ell^1;p),\ldots,\nu_f(\ell^k;p)\}>\nu_f.
$$
\end{Lemma}

\Proof Let us write $h=g_1\ell^1+\cdots+g_k\ell^k$ for suitable
$g_1,\ldots,g_k\in\ca
O_{M,p}$. Then
$$
h\circ f-h=\sum_{j=1}^k\bigl[(g_j\circ f)(\ell^j\circ f-\ell^j)+
    (g_j\circ f-g_j)\ell^j\bigr],
$$
and the assertion follows.\hfill\qed 

\begin{Cor} \label{nondege}  If $f$ is tangential at one point~$p\in S${\rm ,}
then it is tangential at all points of~$S$.
\end{Cor}

\Proof The coherence of the sheaf of ideals of~$S$ implies that if
$\{\ell^1,\ldots,\ell^k\}$
are generators of~$\ca I_{S,p}$ then the corresponding germs are generators
of $\ca I_{S,q}$
for all $q\in S$ close enough to~$p$. Then Lemmas~\ref{fattouno}.(ii)
and~\ref{fattodue}  imply that
both the set of points where $f$ is tangential and the set of points
where $f$ is not tangential are open; hence the assertion follows
because $S$ is connected.\Endproof \vglue4pt

Of course, we shall then say that $f$ is {\it tangential} along~$S$ if
it is tangential at any point of~$S$.

\begin{Example} Let $p$ be a smooth point of~$S$, and
choose local coordinates $z=(z^1,\ldots,z^n)$ defined in a
neighborhood~$U$ of~$p$, centered at~$p$ and such that $S\cap
U=\{z^1=\cdots=z^m=0\}$. We shall write $z'=(z^1,\ldots,z^m)$ and
$z''=(z^{m+1},\ldots,z^n)$, so that $z''$ yields local coordinates
on~$S$. Take $f\in\End(M,S)$, $f\not\equiv\id_M$; then in local
coordinates the map~$f$ can be written as $(f^1,\ldots,f^n)$ with
$$
f^j(z)=z^j+\sum_{h\ge 1}P^j_h(z',z''),
$$
where each $P^j_h$ is a homogeneous polynomial of degree~$h$ in the
variables~$z'$, with
coefficients depending holomorphically on~$z''$. Then Lemma~\ref{fattouno}\
yields
$$
\nu_f=\min\{h\ge 1\mid \exists\, 1\le j\le n: P^j_h\not\equiv 0\}.
$$
Furthermore, $\{z^1,\ldots,z^m\}$ is a set of generators of~$\ca
I_{S,p}$; therefore by Lemma~\ref{fattodue}  the map $f$ is tangential
if and only if
$$
\min\{h\ge 1\mid\exists\,1\le j\le m:
P^j_h\not\equiv 0\}>
    \min\{h\ge 1\mid\exists\,m+1\le j\le n: P^j_h\not\equiv 0\}.
$$
\end{Example}

\begin{remark}  When $S$ is smooth, the differential of $f$ acts
linearly on the normal bundle $N_S$ of~$S$ in~$M$. If $S$ is a 
hypersurface, $N_S$ is a
line bundle, and the action is   multiplication by a holomorphic
function~$b$; if
$S$ is compact, this function is a constant. It is easy  to check that in
local
coordinates chosen as in the previous example the expression of the
function~$b$ is exactly
$1+P^1_1(z)/z^1$; therefore we must have
$P^1_1(z)=(b_f-1)z^1$ for a suitable constant~$b_f\in\C$. In particular, if
$b_f\ne 1$ then
necessarily $\nu_f=1$ and
$f$ is not tangential along~$S$.
\end{remark}

\begin{remark}  The number~$\mu$ introduced in~[BT, (2)] is, by Lemma~\ref{fattouno}, our
order of contact; therefore our notion of tangential is equivalent to
the notion of nondegeneracy defined in~[BT] when $n=2$
and~$m=1$. On the other hand, as already remarked in~[BT], a nondegenerate
map in
the sense defined in~[A2] when $n=2$, $m=1$ and $S$ is smooth is tangential
if and only if 
$b_f=1$ (which was the case mainly considered in that paper).
\end{remark}

\begin{Example} A particularly interesting example (actually,
the one inspiring this paper) of map $f\in\End(M,S)$ is obtained by
blowing up a map tangent to the identity. Let $f_o$ be a (germ of)
holomorphic self-map of~$\C^n$ (or of any complex $n$-manifold)
fixing the origin (or any other point) and {\it tangent to the
identity,} that is, such that $d(f_o)_O=\id$. If $\pi\colon
M\to\C^n$ denotes the blow-up of the origin, let
$S=\pi^{-1}(O)\cong\P^{n-1}(\C)$ be the exceptional divisor, and
$f\in\End(M,S)$ the lifting of~$f_o$, that is, the unique
holomorphic self-map of~$M$ such that $f_o\circ\pi=\pi\circ f$
(see, e.g.,~[A1] for details). If
$$
f_o^j(w)=w^j+\sum_{h\ge 2}Q_h^j(w)
$$
is the expansion of $f_o^j$ in a  series of homogeneous polynomials (for
$j=1,\ldots,n$), then
in the canonical coordinates centered in~$p=[1:0:\cdots:0]$ the map $f$ is
given by
$$
f^j(z)=\left\{\begin{array}{ll} {\displaystyle z^1+\sum_{h\ge 2}Q_h^1(1,z'')(z^1)^h} &\hbox{for
$j=1$},\\
\noalign{\vskip6pt}
 {\displaystyle  z^j+{\sum_{h\ge
2}\left[Q_h^j(1,z'')-z^jQ_h^1(1,z'')\right](z^1)^{h-1}
\over 1+\sum_{h\ge 2}Q^1_h(1,z'')(z^1)^{h-1}}}&\hbox{for $j=2,\ldots,n$,}\end{array}\right.
$$
where $z''=(z^2,\ldots,z^n)$. Therefore $b_f=1$, $$\nu_f(z^1;p)=\min\{h\ge
2\mid
Q^1_h(1,z'')\not\equiv 0\},$$ and
\begin{eqnarray*}
\nu_f&=&\min\bigl\{\nu_f(z^1;p),\\
&&\phantom{\big\{\min}\min\{h\ge 1\mid\exists\,2\le j\le n:
Q^j_{h+1}(1,z'')-z^jQ^1_{h+1}(1,z'')\not\equiv 0\}\bigr\}.
\end{eqnarray*}
Let $\nu(f_o)\ge 2$ be the order of~$f_o$, that is, the minimum $h$ such that
$Q^j_h\not\equiv 0$ for some~$1\le j\le n$. Clearly,
$\nu_f(z^1;p)\ge\nu(f_o)$ and
$\nu_f\ge\nu(f_o)-1$. More precisely, if there is   $2\le j\le n$ such that
$Q^j_{\nu(f_o)}(1,z'')\not\equiv z^jQ^1_{\nu(f_o)}(1,z'')$, then
$\nu_f=\nu(f_o)-1$ and $f$
is tangential. If on the other hand we have $Q^j_{\nu(f_o)}(1,z'')\equiv
z^jQ^1_{\nu(f_o)}(1,z'')$ for all $2\le j\le n$, then necessarily
$Q^1_{\nu(f_o)}(1,z'')\not\equiv 0$, $\nu_f(z^1;p)=\nu(f_o)=\nu_f$, and $f$
is
not tangential.
\end{Example}

Borrowing a term from continuous dynamics, we say that a map $f_o$ tangent
to the identity
at the origin is {\it dicritical} if $w^hQ^k_{\nu(f_o)}(w)\equiv
w^kQ^h_{\nu(f_o)}(w)$ for
all $1\le h, k\le n$. Then we have proved that:

\begin{prop} \label{dint} Let $f_o\in\End(\C^n,O)$ be a {\rm (}\/germ of\/{\rm )}   
holomorphic self\/{\rm -}\/map of~$\C^n$  tangent to the identity at the
origin{\rm ,} and let $f\in\End(M,S)$ be its blow\/{\rm -}\/up. Then $f$ is not
tangential if and only if $f_o$ is dicritical. Furthermore{\rm ,}
$\nu_f=\nu(f_o)-1$ if $f_o$ is not dicritical{\rm ,} and
$\nu_f=\nu(f_o)$ if $f_o$ is dicritical.
\end{prop}

In particular, most maps obtained with this procedure are tangential.

\section{Comfortably embedded submanifolds}

Up to now $S$ was any complex subvariety of the manifold $M$. However, some
of the proofs
in the following sections do not work in this generality; so this section is
devoted to
describe the kind of properties we shall (sometimes) need on $S$.

Let $E'$ and $E''$ be two vector bundles on the same manifold~$S$. We recall
(see, e.g.,
[Ati, \S 1])  that an {\it extension} of $E''$ by~$E'$ is an exact sequence of
vector bundles
$$
O\mapright{} E'\mapright{\iota} E\mapright{\pi} E''\mapright{} O.
$$
Two extensions are {\it equivalent} if there is an isomorphism of exact
sequences which is
the
identity on~$E'$ and~$E''$.

A {\it splitting} of an extension $O\mapright{} E'\mapright{\iota}
E\mapright{\pi} E''\mapright{} O$ is a morphism\break $\sigma\colon
E''\to E$ such that $\pi\circ\sigma=\id_{E''}$.  In particular,
$E=\iota(E')\oplus \sigma(E'')$, and we shall say that the
extension {\it splits.} We explicitly remark that an extension
splits if and only if it is equivalent to the trivial extension
$O\to E'\to E'\oplus E''\to E''\to O$.

Let  $S$ now  be a complex submanifold of a complex manifold~$M$. We shall
denote by $TM|_S$
the restriction to~$S$ of the tangent bundle of~$M$, and by
$N_S=TM|_S/TS$ the normal bundle of~$S$ into~$M$. Furthermore, $\ca T_{M,S}$
will be the
sheaf of
germs of holomorphic sections of~$TM|_S$ (which is different from the
restriction~$\ca
T_M|_S$ to~$S$ of the sheaf of holomorphic sections of~$TM$), and $\ca N_S$
the sheaf of
germs of holomorphic sections of~$N_S$.

\begin{defi}  Let $S$ be a complex submanifold of codimension~$m$ in an
$n$-dimensional complex
manifold~$M$. A chart $(U_\alpha,z_\alpha)$ of~$M$ is {\it adapted} to~$S$ if either $S\cap
U_\alpha=\emptyset$ or
$S\cap U_\alpha=\{z^1_\alpha=\cdots=z^m_\alpha=0\}$, where
$z_\alpha=(z^1_\alpha,\ldots,z^n_\alpha)$. In particular,
$\{z^1_\alpha,\ldots,z^m_\alpha\}$ is a set of generators of $\ca
I_{S,p}$ for all $p\in S\cap U_\alpha$. An atlas $\gt
U=\{(U_\alpha,z_\alpha)\}$ of~$M$ is {\it adapted} to~$S$ if all
charts in~$\gt U$ are. If $\gt U=\{(U_\alpha,z_\alpha)\}$ is
adapted to~$S$ we shall denote by $\gt
U_S=\{(U''_\alpha,z''_\alpha)\}$ the atlas of~$S$ given by
$U''_\alpha=U_\alpha\cap S$ and
$z''_\alpha=(z^{m+1}_\alpha,\ldots,z^n_\alpha)$, where we are
clearly considering only the indices such that $U_\alpha\cap
S\ne\emptyset$. If $(U_\alpha,z_\alpha)$ is a chart adapted to~$S$, we
shall denote by $\de_{\alpha,r}$ the projection of $\de/\de
z^r_\alpha|_{S\cap U_\alpha}$ in~$N_S$, and by~$\omega_\alpha^r$
the local section of~$N_S^*$ induced by~$dz^r_\alpha|_{S\cap
U_\alpha}$; thus $\{\de_{\alpha,1},\ldots,\de_{\alpha,m}\}$ and
$\{\omega^1_\alpha,\ldots,\omega^m_\alpha\}$ are local frames for
$N_S$ and $N_S^*$ respectively over~$U_\alpha\cap S$, dual to each
other.
\end{defi}

From now on, every chart and atlas we consider on $M$ will be
adapted to~$S$.

\begin{remark}  We shall use the  Einstein convention on
the sum over repeated indices. Furthermore, indices like
$j$,~$h$,~$k$ will run from~1 to~$n$; indices like
$r$,~$s$,~$t$,~$u$,~$v$ will run from~$1$ to~$m$; and indices
like~$p$,~$q$ will run from~$m+1$ to~$n$.
\end{remark}

\begin{defi}  We shall say that $S$ {\it splits into~$M$} if the extension $O\to
TS\to TM|_S\to
N_S\to O$ splits.
\end{defi}

\begin{Example} It is well-known that if $S$ is a rational smooth curve
 with negative self-intersection in a surface $M$, then $S$
splits into~$M$.
\end{Example}

\begin{prop} \label{splitting} Let $S$ be a complex submanifold of
codimension~$m$ in an $n$\/{\rm -}\/dimensional complex manifold~$M$. Then
$S$ splits into $M$ if and only if there is an atlas $\hat{\gt
U}=\{(\hat U_\alpha,\hat z_\alpha)\}$ adapted to~$S$ such that
\begin{equation}
\left.{\de \hz^p_\beta\over\de \hz^r_\alpha}\right|_S\equiv 0,
\label{eqdbuno} 
\end{equation}
for all $r=1,\ldots, m$, $p=m+1,\ldots, n$ and indices $\alpha$
and $\beta$.
\end{prop}

\Proof  It is well known (see, e.g., [Ati, Prop.~2]) that there
is a one-to-one correspondence between equivalence classes of
extensions of $N_S$ by $TS$ and the cohomology group
$H^1\bigl(S,\Hom(\ca N_S,\ca T_S)\bigr)$, and an extension splits
if and only if it corresponds to the zero cohomology class.

The class corresponding to the extension $O\to TS\to TM|_S\to
N_S\to O$ is the class $\delta(\id_{N_S})$, where $\delta\colon
H^0\bigl(S,\Hom(\ca N_S,\ca N_S)\bigr)\to H^1\bigl(S,\Hom(\ca
N_S,\ca T_S)\bigr)$ is the connecting homomorphism in the long
exact sequence of cohomology associated to the short exact
sequence obtained by applying  the functor~$\Hom(\ca N_S,\cdot)$
to the extension sequence. More precisely, if $\gt U$ is an atlas
adapted to $S$, we get a local splitting morphism
$\sigma_\alpha\colon N_{U''_\alpha}\to TM|_{U''_\alpha}$ by
setting $\sigma_\alpha(\de_{r,\alpha})=\de/\de z^r_\alpha$, and
then the element of $H^1\bigl(\gt U_S,\Hom(\ca N_S,\ca T_S)\bigr)$
associated to the extension is $\{\sigma_\beta-\sigma_\alpha\}$.
Now,
$$
(\sigma_\beta-\sigma_\alpha)(\de_{r,\alpha})=\left.{\de z^s_\beta\over\de
z^r_\alpha}
\right|_S{\de\over\de z^s_\beta}-
    {\de\over\de z^r_\alpha}={\de z^s_\beta\over\de z^r_\alpha}\left.{\de
z^p_\alpha\over\de
z^s_\beta}
\right|_S{\de\over\de z^p_\alpha}.
$$

So, if  \pref{eqdbuno}  holds, then $S$ splits into $M$. Conversely, assume that
$S$ splits into
$M$; then we can find an atlas~$\gt U$ adapted to~$S$ and a 0-cochain
$\{c_\alpha\}\in
H^0(\gt U_S,\ca T_S\otimes\ca N_S^*)$ such that
\begin{equation}
{\de z^s_\beta\over\de z^r_\alpha}\left.{\de z^p_\alpha\over\de
z^s_\beta}\right|_S
    =(c_\beta)^q_s{\de z^s_\beta\over\de z^r_\alpha}\left.{\de
z^p_\alpha\over\de
z^q_\beta}\right|_S-(c_\alpha)^p_r
\label{eqspdue}
\end{equation}
on $U_\alpha\cap U_\beta\cap S$. We claim that the coordinates
\begin{equation}
\left\{ \begin{array}{l} \hat z^r_\alpha= z^r_\alpha,\\ \noalign{\vskip4pt}
 \hat z^p_\alpha=
z^p_\alpha+(c_\alpha)^p_s(z''_\alpha)
    z^s_\alpha \end{array}\right.
\label{eqsplittingatlas}
\end{equation}
satisfy \pref{eqdbuno}  when restricted to suitable open sets $\hat
U_\alpha\subseteq U_\alpha$.
Indeed, \pref{eqspdue} yields
\begin{eqnarray*}
{\de\hat z^p_\beta\over\de\hat z^r_\alpha}&=&{\de\hat
z^p_\beta\over\de z^s_\alpha}
    {\de z^s_\alpha\over\de\hat z^r_\alpha}+{\de\hat z^p_\beta\over\de
z^q_\alpha}
    {\de z^q_\alpha\over\de\hat z^r_\alpha}={\de\hat z^p_\beta\over\de
z^r_\alpha}-(c_\alpha)^q_r {\de\hat z^p_\beta\over\de z^q_\alpha}+R_1\\ 
&=&{\de z^p_\beta\over\de z^r_\alpha}+(c_\beta)^p_s
{\de z^s_\beta\over\de z^r_\alpha}-(c_\alpha)^q_r
{\de z^p_\beta\over\de z^q_\alpha}+R_1=R_1,
\end{eqnarray*}
where $R_1$ denotes terms vanishing on~$S$, and we are done.\Endproof

\begin{defi}  Assume that $S$ splits into $M$. An atlas $\gt
U=\{(U_\alpha,z_\alpha)\}$
adapted to~$S$ and satisfying~\pref{eqdbuno} \  will be called a {\it splitting atlas} for $S$. It is easy to see that for any splitting
morphism $\sigma\colon N_S\to TM|_S$ there exists a splitting
atlas~$\gt U$ such that $\sigma(\de_{r,\alpha})=\de/\de
z^r_\alpha$ for all $r=1,\ldots m$ and indices~$\alpha$; we shall
say that $\gt U$ is {\it adapted} to~$\sigma$.
\end{defi}

\begin{Example} A {\it local holomorphic retraction} of $M$ onto $S$ is a
holomorphic
retraction $\rho\colon W\to S$, where $W$ is a neighborhood of $S$ in~$M$.
It is clear that the existence of such a local holomorphic retraction
implies that
$S$ splits into~$M$.
\end{Example}

\begin{Example}  Let $\pi\colon M\to S$ be a rank $m$ holomorphic vector
bundle on~$S$.
If we identify $S$ with the zero section of the vector bundle, $\pi$ becomes
a (global)
holomorphic retraction of~$M$ on~$S$. The charts given by the trivialization
of the bundle
clearly give a splitting atlas. Furthermore, if $(U_\alpha,z_\alpha)$ and
$(U_\beta,z_\beta)$ are two such charts, we have
$z''_\beta=\phe_{\beta\alpha}(z''_\alpha)$
and $z'_\beta=a_{\beta\alpha}(z''_\alpha)z'_\alpha$, where $a_{\beta\alpha}$
is an
invertible matrix depending only on $z''_\alpha$. In particular, we have
$$
{\de
z^p_\beta\over\de z^r_\alpha}\equiv 0\qquad\hbox{and}\qquad{\de^2
z^r_\beta\over\de
z^s_\alpha\de z^t_\alpha}\equiv 0
$$
for all $r$,~$s$,~$t=1,\ldots,m$, $p=m+1,\ldots,n$ and indices
$\alpha$ and $\beta$.
\end{Example}

The previous example, compared with \pref{eqdbuno},  suggests the following

\begin{defi}  Let $S$ be a codimension $m$ complex submanifold of an\break
$n$-dimensional complex
manifold~$M$. We say that $S$ is {\it comfortably embedded} in~$M$ if
 $S$ splits into $M$ and
there exists a splitting atlas $\gt U=\{(U_\alpha,z_\alpha)\}$
such that
\begin{equation}
\left.{\de^2 z^r_\beta\over\de z^s_\alpha\de
z^t_\alpha}\right|_S\equiv 0 \label{comfeb}
\end{equation}
for all
$r$,~$s$,~$t=1,\ldots,m$ and indices $\alpha$ and $\beta$.
\end{defi}

 An atlas satisfying the previous condition is said to be  {\it comfortable} for~$S$.
Roughly speaking, then, a comfortably embedded submanifold is like a
first-order
approximation of the zero section of a vector bundle.

Let us express condition \pref{comfeb}  in a different way. If
$(U_\alpha,z_\alpha)$ and $(U_\beta,z_\beta)$ are two charts
about~$p\in S$ adapted to~$S$, we can write
\begin{equation}
z^r_\beta=(a_{\beta\alpha})^r_s\, z^s_\alpha
\label{eqalpha}
\end{equation}
for suitable $(a_{\beta\alpha})^r_s\in\ca O_{M,p}$. The germs
$(a_{\beta\alpha})^r_s$
(unless $m=1$) are not uniquely determined by~\pref{eqalpha}; indeed, all the
other solutions
of~\pref{eqalpha}  are of the form $(a_{\beta\alpha})^r_s+e^r_s$, where the
$e^r_s$'s are
holomorphic and satisfy
\begin{equation}
e^r_s z^s_\alpha\equiv 0.
\label{eqdiffzero}
\end{equation}
Differentiating with respect to $z^t_\alpha$ we get
\begin{equation}
e^r_t+{\de e^r_s\over\de z^t_\alpha}\,z^s_\alpha\equiv 0;
\label{eqdiffuno}
\end{equation}
in particular, $e^r_t|_S\equiv 0$, and so the restriction of
$(a_{\beta\alpha})^r_s$ to
$S$ is uniquely determined --- and it indeed gives the 1-cocycle of the
normal bundle~$N_S$
with respect to the atlas $\gt U_S$. 

Differentiating \pref{eqdiffuno}  we obtain
\begin{equation}
{\de e^r_t\over\de z^s_\alpha}+{\de e^r_s\over\de z^t_\alpha}+{\de^2
e_u^r\over\de
    z^s_\alpha\de z^t_\alpha}\,z^u_\alpha\equiv0;
\label{eqdiffdue}
\end{equation}
in particular,
$$
\left.\left[{\de e^r_t\over\de z^s_\alpha}+{\de e^r_s\over\de
z^t_\alpha}\right]\right|_S
    \equiv 0,
$$
and so the restriction of
$$
{\de(a_{\beta\alpha})^r_t\over\de
z^s_\alpha}+{\de(a_{\beta\alpha})^r_s\over\de z^t_\alpha}
$$
to~$S$ is uniquely determined for all $r$,~$s$,~$t=1,\ldots,m$.

With this notation, we have
$$
{\de^2 z^r_\beta\over\de z^s_\alpha\de
z^t_\alpha}={\de(a_{\beta\alpha})^r_s\over\de
    z^t_\alpha}+{\de(a_{\beta\alpha})^r_t\over\de
z^s_\alpha}+{\de^2(a_{\beta\alpha})^r_u
    \over\de z^s_\alpha\de z^t_\alpha}z^u_\alpha;
$$
therefore \pref{comfeb}  is equivalent to requiring
\begin{equation}
\left.\left({\de(a_{\beta\alpha})^r_t\over\de
z^s_\alpha}+{\de(a_{\beta\alpha})^r_s\over\de
z^t_\alpha}\right)\right|_S\equiv 0
\label{eqdbdue}
\end{equation}
for all $r$,~$s$,~$t=1,\ldots,m$, and indices $\alpha$ and
$\beta$.

\begin{Example}  It is easy to check that the exceptional divisor~$S$ in
Example~1.2 is
comfortably embedded into the blow-up~$M$.
\end{Example}

Then the main result of this section is

\begin{theorem} \label{ce}  Let $S$ be a codimension $m$ complex submanifold of an\break
$n$\/{\rm -}\/dimensional complex manifold~$M$. Assume that $S$ splits into $M${\rm ,} and
let $\gt
U=\{(U_\alpha,z_\alpha)\}$ be a splitting atlas. Define a $1$\/{\rm -}\/cochain
$\{h_{\beta\alpha}\}$
of~$\ca N_S\otimes \ca N_S^*\otimes\ca N_S^*$ by setting
\begin{eqnarray}
h_{\beta\alpha}&=&{1\over2}\left.{\de z^r_\alpha\over\de z^u_\beta}
{\de^2 z^u_\beta\over\de z^s_\alpha\de z^t_\alpha}\right|_S
    \de_{\alpha,r}\otimes \omega_\alpha^s\otimes
\omega_\alpha^t \label{eqh}
 \\ 
&=&{1\over2}(a_{\alpha\beta})^r_u\left.\left(
    {\de(a_{\beta\alpha})^u_s\over\de
z^t_\alpha}+{\de(a_{\beta\alpha})^u_t\over\de
z^s_\alpha}\right)\right|_S\de_{\alpha,r}\otimes \omega_\alpha^s\otimes
\omega_\alpha^t.\nonumber
\end{eqnarray}
Then\/{\rm :}\/
\begin{itemize}
\ritem{\rm (i)} $\{h_{\beta\alpha}\}$ defines an element $[h]\in H^1(S,\ca
N_S\otimes\ca N_S^*\otimes\ca N_S^*)$ independent of~$\gt U$;
\ritem{(ii)} $S$ is comfortably embedded in~$M$ if and only if
$[h]=0$.
\end{itemize}
\end{theorem}

\Proof (i) Let us first prove that $\{h_{\beta\alpha}\}$ is a
1-cocycle with values in $\ca N_S\otimes\ca N_S^*\otimes\ca
N_S^*$. We know that
$$
(a_{\alpha\beta})^r_u(a_{\beta\alpha})^u_s=\delta^r_s+e^r_s,
$$
where $\delta^r_s$ is Kronecker's delta, and the $e^r_s$'s satisfy \pref{eqdiffzero} .
Differentiating we get
$$
{\de(a_{\alpha\beta})^r_u\over\de z^t_\alpha}(a_{\beta\alpha})^u_s+
    (a_{\alpha\beta})^r_u{\de(a_{\beta\alpha})^u_s\over\de z^t_\alpha}
    ={\de e^r_s\over\de z^t_\alpha};
$$
therefore~\pref{eqdiffdue}  yields
$$
(a_{\beta\alpha})^u_s\left.{\de(a_{\alpha\beta})^r_u\over\de
z^t_\alpha}\right|_S+
(a_{\beta\alpha})^u_t\left.{\de(a_{\alpha\beta})^r_u\over\de
z^s_\alpha}\right|_S
=-(a_{\alpha\beta})^r_u\left.\left({\de(a_{\beta\alpha})^u_s\over\de
z^t_\alpha}+
    {\de(a_{\beta\alpha})^u_t\over\de z^s_\alpha}\right)\right|_S.
$$
Hence
\begin{eqnarray*}
h_{\alpha\beta}&=&{1\over2}(a_{\beta\alpha})^r_u\left.\left(
    {\de(a_{\alpha\beta})^u_s\over\de
z^t_\beta}+{\de(a_{\alpha\beta})^u_t\over\de
z^s_\beta}\right)\right|_S\de_{\beta,r}\otimes\omega_\beta^s\otimes
\omega_\beta^t\\ 
&=&{1\over2}(a_{\beta\alpha})^r_u(a_{\alpha\beta})^{r_1}_r(a_{\beta\alpha})^s
_{s_1}
    (a_{\beta\alpha})^t_{t_1}\\
&&\times \left.\left((a_{\alpha\beta})^{t_2}_t
    {\de(a_{\alpha\beta})^u_s\over\de
z^{t_2}_\alpha}+(a_{\alpha\beta})^{s_2}_s
{\de(a_{\alpha\beta})^u_t\over\de z^{s_2}_\alpha}\right)\right|_S
    \de_{\alpha,r_1}\otimes \omega_\alpha^{s_1}\otimes
\omega_\alpha^{t_1}\\ 
&=&{1\over2}\left.\left((a_{\beta\alpha})^s_{s_1}{\de(a_{\alpha\beta})^{r_1}_
s\over\de
z^{t_1}_\alpha}+(a_{\beta\alpha})^t_{t_1}{\de(a_{\alpha\beta})^{r_1}_t\over
\de
z^{s_1}_\alpha}\right)\right|_S\de_{\alpha,r_1}\otimes
\omega_\alpha^{s_1}\otimes \omega_\alpha^{t_1}\\ 
&=&-h_{\beta\alpha},
\end{eqnarray*}
where in the second equality we used \pref{eqdbuno}.
Analogously one proves that
$h_{\alpha\beta}+h_{\beta\gamma}+h_{\gamma\alpha}=0$, and thus
$\{h_{\beta\alpha}\}$ is a 1-cocycle as claimed.

Now we have to prove that the cohomology class~$[h]$ is independent of the
atlas~$\gt U$.
Let $\hat{\gt U}=\{(\hat U_\alpha,\hz_\alpha)\}$ be another splitting atlas;
up to taking
a common refinement we can assume that $U_\alpha=\hat U_\alpha$ for
all~$\alpha$.
Choose $(A_\alpha)^r_s\in\ca O(U_\alpha)$ so that
$\hz^r_\alpha=(A_\alpha)^r_s z^s_\alpha$; as usual, the restrictions to $S$
of
$(A_\alpha)^r_s$ and of
$$
{\de(A_\alpha)^r_s\over\de z^t_\alpha}+{\de(A_\alpha)^r_t\over\de
z^s_\alpha}
$$
are uniquely defined. Set, now,
$$
C_\alpha={1\over 2}(A_\alpha^{-1})^r_u
\left.\left[{\de(A_\alpha)^u_s\over\de z^t_\alpha}
+{\de(A_\alpha)^u_t\over\de z^s_\alpha}\right]\right|_S
\de_{\alpha,r}\otimes
    \omega^s_\alpha\otimes\omega^t_\alpha;
$$
then it is not difficult to check that
$$
h_{\beta\alpha}-\hat h_{\beta\alpha}=C_\beta-C_\alpha,
$$
where $\{\hat h_{\beta\alpha}\}$ is the 1-cocycle built using $\hat{\gt U}$,
and this
means exactly that both $\{h_{\beta\alpha}\}$ and $\{\hat h_{\beta\alpha}\}$
determine
the same cohomology class.

\smallskip
(ii) If $S$ is comfortably embedded, using a comfortable atlas we
immediately see that $[h]=0$. Conversely, assume that $[h]=0$; therefore we
can find a
splitting atlas $\gt U$ and a 0-cochain $\{c_\alpha\}$ of $\ca N_S\otimes\ca
N_S^*\otimes\ca N_S^*$ such that
$h_{\beta\alpha}=c_\alpha-c_\beta$. Writing
$$
c_\alpha=(c_\alpha)^r_{st}\,\de_{\alpha,r}\otimes \omega_\alpha^s\otimes
\omega_\alpha^t,
$$
with $(c_\alpha)^r_{ts}$ symmetric in the lower indices, we define $\hz_\alpha$ by setting
$$
\left\{
\begin{array}{ll}
\hz^r_\alpha=z^r_\alpha+(c_\alpha)^r_{st}(z''_\alpha)\,z_\alpha^s
z_\alpha^t& \hbox{for
    $r=1,\ldots,m$,}\\ 
\hz^p_\alpha=z^p_\alpha&\hbox{for $p=m+1,\ldots,n$,}
\end{array}\right.
$$
on a suitable $\hat U_\alpha\subseteq U_\alpha$. Then $\hat{\gt U}=\{(\hat
U_\alpha,\hz_\alpha)\}$ clearly is a splitting atlas; we claim that it is
comfortable too. Indeed, by definition the functions
$$
(\hat
a_{\beta\alpha})^r_s=[\delta^r_u+(c_\beta)^r_{uv}(a_{\beta\alpha})^v_t
z^t_\alpha]
    (a_{\beta\alpha})^u_{u_1} d^{u_1}_s
$$
satisfy \pref{eqalpha} \ for $\hat{\gt U}$, where the $d^{u_1}_s$'s are such that
$z^{u_1}_\alpha=d^{u_1}_s\hz^s_\alpha$. Hence
\begin{eqnarray*}
\left.\left({\de(\hat a_{\beta\alpha})^r_s\over\de\hz^t_\alpha}+
    {\de(\hat a_{\beta\alpha})^r_t\over\de\hz^s_\alpha}\right)\right|_S&=&
2(c_\beta)^r_{uv}(a_{\beta\alpha})^u_s(a_{\beta\alpha})^v_t|_S+
\left.\left({\de(a_{\beta\alpha})^r_s\over\de z^t_\alpha}+
    {\de(a_{\beta\alpha})^r_t\over\de z^s_\alpha}\right)\right|_S\\
&&
+(a_{\beta\alpha})^r_u\left.\left({\de d^u_s\over\de z^t_\alpha}+{\de
d^u_t\over\de
z^s_\alpha}\right)\right|_S.
\end{eqnarray*}
Now, differentiating
$$
z^u_\alpha=d^u_v\bigl(z^v_\alpha+(c_\alpha)^v_{rs}z^r_\alpha
z^s_\alpha\bigr)
$$
we get
$$
\delta^u_t={\de d^u_v\over\de
z^t_\alpha}\bigl(z^v_\alpha+(c_\alpha)^v_{rs}z^r_\alpha
z^s_\alpha\bigr)+d^u_v\bigl(\delta^v_t+2(c_\alpha)^v_{rt}z^r_\alpha\bigr)
$$
and
$$
0=\left.\left({\de d^u_s\over\de z^t_\alpha}+{\de d^u_t\over\de
z^s_\alpha}\right)\right|_S
+2(c_\alpha)^u_{st}.
$$
Recalling that $h_{\beta\alpha}=c_\alpha-c_\beta$ we then see that $\hat{\gt
U}$ satisfies
\pref{eqdbdue}, and we are done.\Endproof 

\begin{remark}  Since $N_S\otimes N_S^*\otimes
N_S^*\cong\Hom\bigl(N_S,\Hom(N_S,N_S)\bigr)$,
the previous theorem asserts that to any submanifold $S$ splitting into~$M$
we can
canonically associate an extension
$$
O\to\Hom(N_S,N_S)\to E\to N_S\to O
$$
of $N_S$ by $\Hom(N_S,N_S)$, and $S$ is comfortably embedded in
$M$ if and only if this extension splits.  See also [ABT] for more details on comfortably embedded submanifolds.
\end{remark}

\section{The canonical sections}

Our next aim is to associate to any $f\in\End(M,S)$ different from the
identity a
section of a suitable vector bundle, indicating (very roughly speaking) how
$f$ would
move $S$ if it did not keep it fixed. To do so, in this section we still
assume that $S$ is a smooth complex submanifold of a complex manifold~$M$;
however, in Remark~3.3 we shall describe the changes needed to avoid this
assumption.

Given $f\in\End(M,S)$, $f\not\equiv\id_M$, it is clear that $df|_{TS}=\id$;
therefore
$df-\id$ induces a map from $N_S$ to~$TM|_S$, and thus  a
holomorphic section
over~$S$ of the bundle~$TM|_S\otimes N^*_S$.
If $(U,z)$ is a chart adapted to~$S$, we
can define germs $g^h_r$ for $h=1,\ldots,n$ and $r=1,\ldots,m$
by writing
$$
z^h\circ f-z^h=z^1 g^h_1+\cdots+z^m g^h_m.
$$
It is easy to check that the germ of the section of
$TM|_S\otimes N^*_S$ defined by
$df-\id$ is locally expressed by
$$
g^h_r|_{U\cap S}\,
{\de\over\de z_h}\otimes \omega^r,
$$
where we are again indicating by~$\omega^r$ the germ of section of the
conormal bundle
induced by the 1-form $dz^r$ restricted to~$S$.

A problem with this section is that it vanishes identically if (and only if)
$\nu_f>1$. The solution consists in expanding~$f$ at
a higher order.

\begin{defi}  Given a chart $(U,z)$ adapted to~$S$, set
$f^j=z^j\circ f$, and write
\begin{equation}
f^j-z^j= z^{r_1}\cdots z^{r_{\nu_f}}\,g^j_{r_1\ldots r_{\nu_f}},
\label{eqXfprep}
\end{equation}
where the $g^j_{r_1\ldots r_{\nu_f}}$'s are symmetric in $r_1,
\ldots, r_{\nu_f}$ and do not all vanish restricted to~$S$. Let us
then define
\begin{equation}
\ca X_f= g^h_{r_1\ldots r_{\nu_f}}\,
    {\de\over\de z^h}\otimes dz^{r_1}\otimes\cdots\otimes dz^{r_{\nu_f}}.
\label{eqLfdif}
\end{equation}
This is a local section of $TM\otimes(T^*M)^{\otimes\nu_f}$,
defined in a neighborhood of a point of~$S$; furthermore, when
restricted to~$S$, it induces a local section of
$TM|_S\otimes(N_S^*)^{\otimes\nu_f}$.
\end{defi}

\begin{remark}  When $m>1$ the $g^j_{r_1\ldots r_{\nu_f}}$'s are {\it not} uniquely
determined by \pref{eqXfprep}. Indeed, if $e^j_{r_1\ldots r_{\nu_f}}$ are such
that
\begin{equation}
e^j_{r_1\ldots r_{\nu_f}} z^1\cdots z^{r_{\nu_f}}\equiv 0
\label{eqe}
\end{equation}
then $g^j_{r_1\ldots r_{\nu_f}}+e^j_{r_1\ldots r_{\nu_f}}$ still
satisfies \pref{eqXfprep}. This means that the section \pref{eqLfdif} is not
uniquely determined too; but, as we shall see, this will not be a
problem. For instance, \pref{eqe} implies that $e^j_{r_1\ldots
r_{\nu_f}}\in\ca I_S$; therefore $\ca X_f|_{U\cap S}$ is always
uniquely determined --- though {\it a priori} it might depend on the
chosen chart. On the other hand, when $m=1$ both the $g^j_{r_1\ldots
r_{\nu_f}}$'s and $\ca X_f$ {\it are} uniquely determined; this is one of
the reasons
making the codimension-one case simpler than the general case.
\end{remark}

We have already remarked that when $\nu_f=1$ the section
$\ca X_f$ restricted to $U\cap S$ coincides with the restriction of
$df-\id$ to $S$. Therefore when $\nu_f=1$ the restriction of $\ca X_f$
to~$S$ gives a
globally well-defined section. Actually, this holds for any~$\nu_f\ge 1$:

\begin{prop} \label{norm} Let $f\in\End(M,S)$, $f\not\equiv\id_M$. Then the
restriction
of~$\ca X_f$ to~$S$ induces a global holomorphic section~$X_f$ of
the bundle $TM|_S\otimes(N_S^*)^{\otimes\nu_f}$.
\end{prop}

\Proof Let $(U,z)$ and $(\hat U,\hz)$ be two charts about $p\in S$
adapted to $S$. Then we can find holomorphic functions
$a^r_s$ such that
\begin{equation}
\hz^r=a_s^r\, z^s;
\label{eqcambiocoord}
\end{equation}
in particular,
\begin{equation}
{\de \hz^r\over\de z^s}=a^r_s\quad\hbox{(mod $\ca
I_S$)}\qquad\hbox{and}\qquad {\de
\hz^r\over\de z^p}=0\quad\hbox{(mod $\ca I_S$).}
\label{eqabcc}
\end{equation}
Now set $f^j=z^j\circ f$, $\hat f^j=\hz^j\circ f$, and define
$g^j_{r_1\cdots r_{\nu_f}}$
and $\hg^j_{r_1\cdots r_{\nu_f}}$ using
\pref{eqXfprep}  with $(U,z)$ and~$(\hat U,\hz)$ respectively. Then
\pref{eqcambiocoord}  and \pref{eqbasedue}  yield
\begin{eqnarray*}
a_{s_1}^{r_1}\cdots a_{s_{\nu_f}}^{r_{\nu_f}}\,\hg^j_{r_1\ldots
r_{\nu_f}}z^{s_1}
    \cdots z^{s_{\nu_f}}&=&\hg^j_{r_1\ldots
r_{\nu_f}}\hz^{r_1}\cdots\hz^{r_{\nu_f}}\\ 
&=&\hat f^j-\hz^j=(f^h-z^h){\de\hz^j\over\de z^h}+R_{2\nu_f}\\
&=&g^h_{s_1\ldots
s_{\nu_f}}
{\de\hz^j\over\de z^h}z^{s_1}\cdots z^{s_{\nu_f}}+R_{2\nu_f},
\end{eqnarray*}
where the remainder terms $R_{2\nu_f}$ belong to~$\ca I_S^{2\nu_f}$.
Therefore we find
\begin{equation}
a^{r_1}_{s_1}\cdots
a^{r_{\nu_f}}_{s_{\nu_f}}\,\hg^j_{r_1\ldots r_{\nu_f}}=
    {\de\hz^j\over\de z^h}\,g^h_{s_1\ldots s_{\nu_f}}\quad\hbox{(mod $\ca
I_S$).}
\label{eqlocalchange}
\end{equation}
Recalling \pref{eqabcc}  we then get
\begin{eqnarray*}&& \hskip-24pt
\hg^j_{r_1\ldots r_{\nu_f}}\,{\de\over\de\hz^j}\otimes
d\hz^{r_1}\otimes\cdots\otimes d\hz^{r_{\nu_f}} \\
&&\qquad =
{\de z^h\over\de\hz^j}{\de\hz^{r_1}\over\de
z^{k_1}}\cdots{\de\hz^{r_{\nu_f}}\over\de
z^{k_{\nu_f}}}\hg^j_{r_1\ldots r_{\nu_f}}
    \,{\de\over\de z^h}\otimes dz^{k_1}\otimes\cdots\otimes
    dz^{k_{\nu_f}}\\ 
&&\qquad =a^{r_1}_{s_1}\cdots a^{r_{\nu_f}}_{s_{\nu_f}}\hg^j_{r_1\ldots
r_{\nu_f}}{\de
z^h\over\de\hz^j}\,{\de\over\de z^h}\otimes dz^{s_1}\otimes\cdots\otimes
    dz^{s_{\nu_f}}
\quad\hbox{(mod $\ca I_S$)}\\ 
&&\qquad =g^h_{s_1\ldots s_{\nu_f}}{\de\over\de z^h}\otimes
dz^{s_1}\otimes\cdots\otimes
    dz^{s_{\nu_f}}\quad\hbox{(mod $\ca I_S$),}
\end{eqnarray*}
and we are done.\Endproof 

\begin{remark}  For later use, we explicitly notice that when
$m=1$ the germs $a^r_s$ are uniquely determined, and
\pref{eqlocalchange} becomes
\begin{equation}
(a^1_1)^{\nu_f}\hat g^j_{1\ldots 1}={\de\hz^j\over\de
z^h}\,g^h_{1\ldots1} \quad\hbox{(mod $\ca I_S^{\nu_f}$).}
\label{eqlocalchangecoduno}
\end{equation}
\end{remark}

\begin{defi}  Let $f\in\End(M,S)$, $f\not\equiv\id_M$. The {\it canonical section
 $X_f\in
H^0\bigl(S,\,\ca T_{M,S}\otimes(\ca N^*_S)^{\otimes\nu_f}\bigr)$ associated
to~$f$} is
defined by setting
\begin{equation}
X_f=g^h_{s_1\ldots s_{\nu_f}}|_S\,{\de\over\de z^h}\otimes
\omega^{s_1}\otimes\cdots\otimes
    \omega^{s_{\nu_f}}
\label{eqXfdif}
\end{equation}
in any chart adapted to~$S$. Since
$(N_S^*)^{\otimes\nu_f}=(N_S^{\otimes\nu_f})^*$, we can also
think of $X_f$ as a holomorphic section of $\Hom(N_S^{\otimes\nu_f},TM|_S)$,
and introduce the {\it canonical distribution}
$\Xi_f=X_f(N_S^{\otimes\nu_f})\subseteq TM|_S$.
\end{defi}

In particular we can now justify the term ``tangential'' previously
introduced:

\begin{Cor} \label{deftang} Let $f\in\End(M,S)${\rm ,} $f\not\equiv\id_M$. Then
$f$ is
tangential if and only if the canonical distribution is tangent to~$S${\rm ,} that
is if and only if $\Xi_f\subseteq TS$.
\end{Cor}

\Proof  This follows from Lemma~\ref{fattodue}.\Endproof  

\begin{Example} By the notation  introduced in
Example~1.2, if $f$ is obtained by blowing up a map~$f_o$ tangent
to the identity, then the canonical coordinates centered in
$p=[1:0:\cdots:0]$ are adapted to~$S$. In particular, if $f_o$ is
non-dicritical (that is, if $f$ is tangential) then in a
neighborhood of $p$,
$$
X_f=\bigl[Q^q_{\nu(f_o)}(1,z'')-z^qQ^1_{\nu(f_o)}(1,z'')\bigr]\,{\de\over\de
z^q}\otimes (\omega^1)^{\otimes(\nu(f_o)-1)}.
$$
\end{Example}

\begin{remark}  To be more precise, $X_f$ is a section of the subsheaf  ${\cal T}_{M,S}\otimes \hbox{\rm
Sym}^{\nu_f}({\cal N}_S^\ast)$, where ${\rm Sym}^{\nu_f}({\cal N}_S^\ast)$ is the symmetric $\nu_f$-fold
tensor product of ${\cal N}_S^\ast$.  Now, the sheaf
 $\ca N_S^*$ is isomorphic to $\ca
I_S/\ca I_S^2$, and it is known that
${\rm Sym}^{\nu_f}{\cal I}_S/{\cal I}_S^2$
 is isomorphic to ${\ca I}_S^{\nu_f}/{\ca
I}_S^{\nu_f+1}$. This
allows us
to define $X_f$ as a global section of the coherent sheaf $\ca
T_{M,S}\otimes {\rm Sym}^{\nu_f}({\cal I}_S/{\cal I}_S^2)$ even when $S$ is singular. Indeed, if $(U,z)$ is a
local chart
adapted
to~$S$, for $j=1,\ldots,n$ the functions $f^j-z^j$ determine local sections
$[f^j-z^j]$ of
$\ca
I_S^{\nu_f}/\ca I_S^{\nu_f+1}$. But, since for any other chart $(\hat U,\hat
z)$,
$$
 \hat f^j-\hat z^j=(f^h-z^h){\de\hat z^j\over\de z^h}+R_{2\nu_f},
$$
then
$(\de/\de z^j)\otimes[f^j-z^j]$ is a well-defined global section of~$\ca
T_{M,S}\otimes{\rm Sym}^{\nu_f}({\cal I}_S/{\cal I}_S^2)$ which coincides with~$X_f$ when $S$ is smooth.
\end{remark}

\begin{remark}  When $f$ is tangential and $\Xi_f$ is involutive as a
sub-distribution
of~$TS$ --- for instance when $m=1$ --- we thus get a holomorphic singular
foliation on~$S$
canonically associated to~$f$. As already remarked in~[Br], this possibly is
the 
reason explaining the similarities discovered in~[A2] between the local
dynamics of
holomorphic maps tangent to the identity and the dynamics of singular
holomorphic
foliations.
\end{remark}

\begin{defi}  A point $p\in S$ is {\it singular} for $f$ if there exists
$v\in(N_S)_p$, $v\ne O$, such that $X_f(v\otimes\cdots\otimes v)=O$. We
shall denote
by~$\Sing(f)$
the set of singular points of~$f$.
\end{defi}

In Section~7 it will become clear why we choose this definition for
singular points. In Section~8 we shall describe a dynamical
interpretation of~$X_f$ at nonsingular points in the codimension-one case; see Proposition~8.1.

\begin{remark}  If $S$ is a  hypersurface, the normal
bundle is a line bundle. Therefore $\Xi_f$ is a 1-dimensional
distribution, and the singular points of $f$ are the points where
$\Xi_f$ vanishes. Recalling \pref{eqXfdif}, we then see that
$p\in\Sing(f)$ if and only if $g^1_{1\ldots1}(p)=\cdots=g^n_{1\ldots1}(p)=0$
for any adapted chart, and thus both the
strictly fixed points of~[A2] and the singular points of~[BT], [Br]
are singular in our case as well.
\end{remark}

As we shall see later on, our index theorems will need a section of
$TS\otimes(N_S^*)^{\otimes\nu_f}$; so it will be natural to assume $f$
tangential. When $f$
is not tangential but $S$ splits in~$M$ we can work too.

Let $O\mapright{} TS\mapright{\iota} TM|_S\mapright{\pi}
N_S\mapright{} O$ be the usual extension. Then we can associate to
any splitting morphism $\sigma\colon N_S\to TM|_S$ a morphism
$\sigma'\colon TM|_S\to TS$ such that
$\sigma'\circ\iota=\id_{TS}$, by
$\sigma'=\iota^{-1}\circ(\sigma\circ\pi-\id_{TM|_S})$. Conversely,
if there is a morphism $\sigma'\colon TM|_S\to TS$ such that
$\sigma'\circ\iota=\id_{TS}$, we get a splitting morphism by
setting $\sigma=(\pi|_{\Ker\sigma'})^{-1}$. Then

\begin{defi}  Let $f\in\End(M,S)$, $f\not\equiv\id_M$, and assume that $S$ splits
in~$M$. Choose
a splitting morphism $\sigma\colon N_S\to TM|_S$
 and let
$\sigma'\colon TM|_S\to TS$ be the induced morphism. We
 set
$$
H_{\sigma,f}=(\sigma'\otimes\id)\circ X_f\in H^0\bigl(S,\ca
T_S\otimes(N_S^*)^{\otimes\nu_f}\bigr).
$$
Since the differential of $f$ induces a morphism from $N_S$
into itself, we have a dual morphism~$(df)^*\colon N_S^*\to N_S^*$. Then if
$\nu_f=1$
we  also set
$$
H^1_{\sigma,f}=\bigl(\id\otimes(df)^*\bigr)\circ H_{\sigma,f}\in
H^0\bigl(S,\ca
T_S\otimes N_S^*\bigr).
$$
\end{defi}

\begin{remark}  We defined $H^1_{\sigma,f}$ only for $\nu_f=1$ because when
$\nu_f>1$ one has $(df)^*=\id$. On the other hand, when $\nu_f=1$
one has $(df)^*=\id$ if and only if $f$ is tangential. Finally, we
have $X_f\equiv H_{\sigma,f}$ if and only if $f$ is tangential,
and $H_{\sigma,f}\equiv O$ if and only if
$\Xi_f\subseteq\Im\sigma=\Ker\sigma'$.
\end{remark}

Finally, if $(U,z)$ is a chart in an atlas adapted to the splitting
$\sigma$, locally we
have
$$
H_{\sigma,f}=g^p_{s_1\ldots s_{\nu_f}}|_S\,{\de\over\de z^p}\otimes
\omega^{s_1}\otimes\cdots
\otimes\omega^{s_{\nu_f}},
$$
and, if $\nu_f=1$,
$$
H^1_{\sigma,f}=(\delta^s_r+g^s_r)g^p_s|_S\,{\de\over\de z^p}\otimes
\omega^r.
$$

\section{Local extensions}

As we have already remarked, while $X_f$ is well-defined, its
extension $\ca X_f$ in general is not. However, we shall now derive formulas
showing
how to control the ambiguities in the definition of $\ca X_f$, at least in
the cases that interest us most.

In this section we assume $m=1$, i.e., that $S$ has
codimension one in $M$. To simplify notation  we shall write $g^j$ for
$g^j_{1\ldots 1}$
and $a$ for $a^1_1$. We shall also use the following notation:
\begin{itemize}
\ritem{$\bullet$} $T_1$ will denote any sum of terms of the form $g\,
{\de\over\de
z^p}\otimes dz^{h_1}\otimes\cdots\otimes dz^{h_{\nu_f}}$ with $g\in\ca I_S$;
\ritem{$\bullet$} $R_k$ will denote any local section with coefficients in
$\ca
I_S^k$.
\end{itemize}

\noindent  For instance, if $(U,z)$ and $(\hat U,\hat z)$ are two charts
adapted to~$S$,
\begin{eqnarray}
{\de\over\de\hz^h}\otimes (d\hz^1)^{\otimes\nu_f}
&=&a^{\nu_f}
    {\de z^k\over\de\hz^h}\,{\de\over\de z^k}\otimes
(dz^{1})^{\otimes\nu_f}\label{eqgencomp}\\ 
&& +{\de z^1\over\de\hz^h}a^{\nu_f-1}z^1
    \sum_{\ell=1}^{\nu_f}{\de a\over\de z^{j_\ell}}\,{\de\over\de
z^1}\otimes dz^1 \otimes
    \cdots\nonumber\\
&&\cdots \otimes dz^{j_\ell}\otimes\cdots\otimes dz^{1}+T_1+R_2,\nonumber\end{eqnarray}
where
$$
T_1={\de z^p\over\de\hz^h}a^{\nu_f-1}z^1\sum_{\ell=1}^{\nu_f}{\de a\over\de
z^{j_\ell}}
    \,{\de\over\de z^p}\otimes dz^1 \otimes
    \cdots\otimes dz^{j_\ell}\otimes\cdots\otimes dz^1.
$$

Assume now that $f$ is tangential, and let $(U,z)$ be a chart adapted
to~$S$. We know that $f^1-z^1\in\ca I_S^{\nu_f+1}$, and thus we
can write
$$
f^1-z^1=h^1(z^{1})^{\nu_f+1},
$$
where $h^1$ is uniquely determined. Now, if $(\hat U,\hz)$ is
another chart adapted to~$S$ then
\begin{eqnarray*}
a^{\nu_f+1}\hat h^1
    (z^{1})^{\nu_f+1}&=&\hat f^1-\hz^1=(a\circ f)f^1-a z^1\\ 
&=&a(f^1-z^1)+(a\circ f-a)z^1+(a\circ f-a)(f^1-z^1)\\ 
&=&a(f^1-z^1)+{\de a\over\de z^p}(f^p-z^p)z^1+R_{\nu_f+2}\\ 
&=&\left[a h^1+{\de a\over\de z^p} g^p\right]
    (z^{1})^{\nu_f+1}+R_{\nu_f+2}.
\end{eqnarray*}
Therefore
\begin{equation}
a^{\nu_f+1}\hat h^1=
    a h^1+
    {\de a\over\de z^p}
    g^p+R_1.
\label{eqtanggen}
\end{equation}
Since $g^1=h^1 z^{1}$ we then get
\begin{equation}
a^{\nu_f}\hat g^1=
    a g^1+
    {\de a\over\de z^p}
    g^p z^{1}+R_2,
\label{eqtanggeng}
\end{equation}
which generalizes \pref{eqlocalchange}  when $f$ is tangential and $m=1$.

Putting \pref{eqtanggeng}, \pref{eqlocalchange}  and \pref{eqgencomp} into
\pref{eqLfdif} we then get

\begin{Lemma}\label{ccin} Let $f\in\End(M,S)${\rm ,} $f\not\equiv\id_M$. Assume that
$f$ is
tangential{\rm ,} and that $S$ has codimension~$1$. Let $(\hat U,\hz)$ and $(U,z)$
be two charts
about $p\in S$ adapted to~$S${\rm ,} and let $\hat{\ca X}_f$, $\ca X_f$ be given
by~{\rm \pref{eqLfdif}} in the respective coordinates. Then
$$
\hat{\ca X}_f=\ca X_f+T_1+R_2.
$$
\end{Lemma}

When $S$ is comfortably embedded in $M$ and of codimension one
 we shall also need nice local extensions of
$H_{\sigma,f}$ and $H^1_{\sigma,f}$, and to know how they behave under
change of
(comfortable) coordinates.

\begin{defi}  Let $S$ be comfortably embedded in $M$ and of codimension $1$, and
take $f\in\End(M,S)$, $f\not\equiv\id_M$. Let $(U,z)$ be a chart
in a comfortable atlas, and set $b^1(z)=g^1(O,z'')$; notice that
$f$ is tangential if and only if $b^1\equiv O$. Write
$g^1=b^1+h^1z^1$ for a well-defined holomorphic function~$h^1$;
then set
\begin{equation}
\ca H_{\sigma,f}=h^1z^{1}\,{\de\over\de
z^1}\otimes(dz^1)^{\otimes \nu_f}+ g^p\,{\de\over\de z^p}\otimes
(dz^{1})^{\otimes\nu_f},\label{hestende}
\end{equation}
and if $\nu_f=1$ set
\begin{equation}
\ca H^1_{\sigma,f}=h^1z^{1}\,{\de\over\de z^1}\otimes dz^1+
g^p(1+b^1)\,{\de\over\de z^p}\otimes dz^1.\label{hestendeuno}
\end{equation}
Notice that $\ca H_{\sigma, f}$ (respectively, $\ca
H^1_{\sigma,f}$) restricted to $S$ yields
$H_{\sigma, f}$ (respectively, $ H^1_{\sigma, f}$).
\end{defi}

\begin{prop}  Let $f\in\End(M,S)$, $f\not\equiv\id_M$. Assume
that $S$ is comfortably embedded in $M${\rm ,} and of codimension one.
Fix a comfortable atlas~$\gt U${\rm ,} and let $(U,z)${\rm ,} $(\hat U,\hat
z)$ be two charts in~$\gt U$ about $p\in S$. Then  if~$\nu_f=1${\rm ,}
\begin{equation}
\hat{\ca H}^1_{\sigma,f}=\ca H^1_{\sigma,f}+T_1+R_2,
\end{equation}
 while if $\nu_f>1${\rm ,}
\begin{equation}
\hat{\ca H}_{\sigma, f}=\ca H_{\sigma, f}+T_1+R_2,
\end{equation}
where $T_1=T_1^o+T_1^1$ with
\begin{eqnarray*}
T_1^o&=&{1\over a}g^q z^1\sum_{\ell=1}^{\nu_f}{\de a\over\de
z^{p_\ell}}\,{\de\over\de
z^q}\otimes
dz^1\otimes\cdots\otimes dz^{p_\ell}\otimes\cdots\otimes dz^1,\\
T_1^1&=&-ag^1{\de
z^q\over\de
\hz^1}\, {\de\over\de z^q}\otimes(dz^1)^{\otimes\nu_f}.
\end{eqnarray*}
\end{prop}

\Proof  First of all, from \pref{eqlocalchangecoduno}, $a^{\nu_f}\hat
b^1=a b^1 \hbox{(mod $\ca I_S$)}$. But since we are using a comfortable
atlas we get
$$
{\de (a^{\nu_f}\hat b^1 - ab^1)\over \de z^1}=(\nu_f a^{\nu_f-1}\hat
b^1-b^1)
{\de a\over\de z^1}+R_1=R_1,
$$
and thus
\begin{equation}
a^{\nu_f}\hat
b^1=a b^1 \quad\hbox{(mod $\ca I_S^{2}$)}.
\end{equation}
If $\nu_f>1$ then by \pref{eqlocalchangecoduno} and
(4.8),
$$
a^{\nu_f}\hat{h}^1\hat{z}^1=(ah^1+{\de a\over \de z^p}g^p)z^1
\quad\hbox{(mod $\ca I_S^2$),}
$$
which implies
\begin{equation}
a^{\nu_f+1}\hat{h}^1=ah^1+{\de a\over \de z^p}g^p \quad\hbox{(mod
$\ca I_S$).}
\end{equation}
If $\nu_f=1$, using \pref{comfeb}  we can write
\begin{eqnarray*}
\hat b^1\hz^1+\hat h^1(\hz^{1})^2&=&\hat f^1-\hz^1\\ 
    &=&{\de\hz^1\over\de z^j}(f^j-z^j)+{1\over2}{\de^2\hz^1\over\de z^h\de
z^k}(f^h-z^h)
    (f^k-z^k)+R_3\\ 
&=&ab^1z^1+\left[ah^1+{\de a\over\de
z^p}g^p(1+b^1)\right](z^1)^2+R_3, 
\end{eqnarray*}
and by (4.8),
\begin{equation}
a^2\hat h^1=ah^1+{\de a\over\de z^p}g^p(1+b^1)\quad\hbox{(mod $\ca
I_S$)}.
\label{cambiohuno}
\end{equation}
So if we compute $\hat{\ca H}_{\sigma,f}$ for $\nu_f>1$ (respectively,
$\hat{\ca
H}^1_{\sigma,f}$ for $\nu_f=1$) using \pref{eqlocalchangecoduno}, \pref{eqgencomp} and
(4.9)
(respectively, \pref{eqlocalchangecoduno}, \pref{eqgencomp}, (4.8) and
(4.10)), we get the assertions.\Endproof 

\section{Holomorphic actions}

The index theorems  to be discussed depend on actions of
vector bundles. This concept was introduced by Baum and Bott
in [BB], and later generalized in [CL], [LS], [LS2] and [Su]. Let
us recall here the relevant definitions.

Let  $S$ again be a submanifold of codimension~$m$ in an
$n$-dimensional complex manifold~$M$, and let $\pi_F\colon F\to S$
be a holomorphic vector bundle on~$S$. We shall denote by $\ca F$
the sheaf of germs of holomorphic sections of~$F$, by $\ca T_S$
the sheaf of germs of holomorphic sections of~$TS$, and by
$\Omega^1_S$ (respectively, $\Omega^1_M$) the sheaf of holomorphic 1-forms
on~$S$
(respectively, on $M$).

A section $X$ of $\ca T_S\otimes\ca F^*$ (or, equivalently, a holomorphic
section of\break
$TS\otimes F^*$) can be interpreted as a morphism $X\colon\ca F\to\ca T_S$.
Therefore it induces a derivation $X^\#\colon\ca O_S\to\ca F^*$ by setting
\begin{equation}
X^\#(g)(u)=X(u)(g)
\end{equation}
for any $p\in S$, $g\in\ca O_{S,p}$ and $u\in\ca F_p$. If
$\{f_1^*,\ldots,f_k^*\}$ is a
local frame for~$F^*$ about~$p$, and $X$ is locally given by $X=\sum_j
v_j\otimes f_j^*$,
then
\begin{equation}
X^\#(g)=\sum_jv_j(g)f_j^*.
\end{equation}
Notice that if $X^*\colon\Omega^1_S\to\ca F^*$ denotes the dual morphism of
$X\colon\ca
F\to\ca T_S$, by definition we have
$$
X^*(\omega)(u)=\omega\bigl(X(u)\bigr)
$$
for any $p\in S$, $\omega\in(\Omega^1_S)_p$ and $u\in\ca F_p$, and so
$$
X^\#(g)=X^*(dg).
$$

\begin{defi}  Let $\pi_E\colon E\to S$ be another
holomorphic vector bundle on~$S$, and denote by $\ca E$ the sheaf
of germs of holomorphic sections of~$E$. Let $X$ be a section of
$\ca T_S\otimes\ca F^*$. A {\it holomorphic
action of~$F$ on $E$ along~$X$} (or an {\it $X$-connection} on~$E$) is a
${\Bbb C}$-linear map  $\tilde
X\colon\ca E\to\ca F^*\otimes \ca E$ such that
\begin{equation}
\tilde X(gs)=X^\#(g)\otimes s+g\tilde X(s)
\label{eqaction}
\end{equation}
for any $g\in\ca O_S$ and $s\in\ca E$.
\end{defi}

\begin{Example} If $F=TS$, and the
section~$X$ is the identity $\id\colon TS\to TS$, then $X^\#(g)=dg$, and a
holomorphic action of $TS$ on~$E$ along~$X$ is just a (1,0)-connection
on~$E$.
\end{Example}

\begin{defi}  A point $p\in S$ is a {\it singularity} of a holomorphic section $X$
of $\ca
T_S\otimes\ca F^*$ if the induced map $X_p\colon F_p\to T_pS$ is not
injective. The set
of singular points of~$X$ will be denoted by~$\Sing(X)$, and we shall set
$S^0=S\setminus\Sing(X)$ and $\Xi_X=X(F|_{S^0})\subseteq TS^0$. Notice that
$\Xi_X$ is a
holomorphic subbundle of~$TS^0$.
\end{defi}

The canonical section previously introduced suggests the following
definition:

\begin{defi}  A {\it Camacho-Sad action} on~$S$ is a holomorphic action
of~$N_S^{\otimes\nu}$ on~$N_S$ along a section $X$ of $\ca
T_S\otimes (N_S^{\otimes\nu})^*$, for a suitable $\nu\ge 1$.
\end{defi}

\begin{remark}  The rationale behind the name is the following: as we shall see, the
index
theorem in~[A2] is induced by a holomorphic action of~$N_S^{\otimes\nu_f}$
on~$N_S$
along~$X_f$ when $f$ is tangential, and this index theorem was inspired by the 
Camacho-Sad index theorem [CS].
\end{remark}

Let us describe a way to get Camacho-Sad actions.
Let $\pi\colon TM|_S\to N_S$ be the
canonical projection; we shall use the same symbol for all other projections
naturally induced by it. Let $X$ be any global section
of~$TS\otimes(N_S^{\otimes\nu})^*$. Then we
might try to define an action $\tilde X\colon\ca N_S\to(\ca
N_S^{\otimes\nu})^*\otimes\ca
N_S=\Hom(\ca N_S^{\otimes\nu},\ca N_S)$ by setting
\begin{equation}
\tilde X(s)(u) =\pi([\ca X(\tilde u),\tilde s]|_S)
\end{equation}
for any $s\in\ca N_S$ and $u\in\ca N_S^{\otimes\nu}$, where:
$\tilde s$ is any element in $\ca T_M|_S$ such
that $\pi(\tilde s|_S)=s$; $\tilde u$ is any
element in~$\ca T_M|_S^{\otimes\nu_f}$ such that $\pi(\tilde
u|_S)=u$; and $\ca X$ is a {\it suitably chosen} local section of $\ca
T_M\otimes(\Omega^1_M)^{\otimes\nu}$ that restricted to~$S$ induces~$X$.

Surprisingly enough, we can make this definition work in the cases
interesting to us:

\begin{theorem} \label{duno} Let $f\in\End(M,S)${\rm ,} $f\not\equiv\id_M${\rm ,} be given.
Assume
that $S$ has codimension one in $M$ and that
\vglue4pt
{\rm (a)}  $f$ is tangential to $S${\rm ,} or that
\vglue2pt
{\rm (b)}
$S$  is comfortably embedded into $M$.
\vglue4pt
\noindent Then we can use {\rm (5.4)} to define a Camacho-Sad
action on~$S$ along $X_f$ in case {\rm(a),} along $H_{\sigma,f}$ in case
{\rm(b)}
when~$\nu_f>1${\rm ,} and along $H^1_{\sigma,f}$ in case {\rm(b)} when  $\nu_f=1$.
\end{theorem}

\Proof We shall denote by $X$ the section $X_f$, $H_{\sigma,f}$ or
$H^1_{\sigma,f}$ depending on the case we are considering. Let~$\gt U$
be an atlas adapted to~$S$,
comfortable and adapted to the splitting morphism~$\sigma$ in case (b), and
let $\ca X$ be
the local
extension of~$X$ defined in a chart belonging to~$\gt U$ by Definition~3.1
(respectively,
Definition~4.1). We first prove that the right-hand side of  (5.4) does
not depend on
the chart
chosen. Take
$(U,z)$, $(\hat U,\hz)\in\gt U$ to be  local charts about $p\in S$. Using
Lemma~\ref{ccin} and Proposition~4.2 we get
$$
[\hat {\ca X}(\tilde u),\tilde s]=[(\ca
X+T_1+R_2)(\tilde u),\tilde s] =[\ca{X}(\tilde
u)+T_1+R_2,\tilde s]=[\ca{X}(\tilde u),\tilde s]+T_0+R_1,
$$
where  $T_0$ represents a local section of $TM$ that restricted to~$S$ is
tangent to it. Thus
$$
\pi\bigl([\hat {\ca X}(\tilde u),\tilde
s]|_S\bigr)=\pi\bigl([\ca{X}(\tilde u),\tilde s]|_S\bigr),
$$
as desired.

We must now  show that the right-hand side of  (5.4) does not depend on
the
extensions of $s$ and $u$ chosen. If $\tilde s'$ and $\tilde
u'$ are other extensions of $s$ and $u$ respectively, we have
$(\tilde s'-\tilde s)|_S=T_0$, while $(\tilde u'-\tilde u)|_S$ is a sum of
terms of the
form $V_1\otimes\cdots\otimes V_{\nu_f}$ with at least one $V_\ell$ tangent
to~$S$.
Therefore $\ca X(\tilde u'-\tilde u)|_S=O$ and
\begin{eqnarray*}
[\ca{X}(\tilde u'),\tilde s']|_S&=&[\ca{X}(\tilde
u),\tilde s]|_S+[\ca{X}(\tilde u),\tilde s'-\tilde
s]|_S+[\ca{X}(\tilde u'-\tilde u),\tilde
s]|_S\\
&& +[\ca{X}(\tilde u'-\tilde u),\tilde s'-\tilde
s]|_S=[\ca{X}(\tilde u),\tilde s]|_S+T_0,
\end{eqnarray*}
so that $\pi\bigl([\ca X(\tilde
u'),\tilde s']|_S\bigr)=\pi\bigl([\ca X_f(\tilde u),\tilde s]|_S\bigr)$, as
wanted.

We are left to show that $\tilde X$ is actually an action. Take
$g\in\ca O_S$, and let $\tilde g\in\ca O_M|_S$ be any extension. First of
all,
$$
\tilde X(s)(gu)=\pi\bigl([\ca X(\tilde g\tilde u),\tilde
s]|_S\bigr)=g\tilde X(s)(u)-\tilde s(\tilde
g)|_S\pi\bigl(X(u)\bigr)=g\tilde X(s)(u),
$$
and so $\tilde X(s)$ is a morphism. Finally,
(5.1) yields
$$
\ca X(\tilde u)(\tilde g)|_S=X^\#(g)(u),
$$
and so
$$
\tilde X(gs)(u)=\pi\bigl([\ca X(\tilde u),\tilde g\tilde
s]|_S\bigr)=g\tilde X(s)(u)+
    \ca X(\tilde u)(\tilde g)|_S\,s=g\tilde X(s)(u)+X^\#(g)(u)s,
$$
and we are done. \hfill\qed 

\begin{remark}  If $\nu_f=1$ and $f$ is not tangential then  (5.4) with $\ca
X=\ca H_{\sigma,f}$
does not define an action. This is the reason why we introduced the
new section $H^1_{\sigma,f}$ and its extension $\ca
H^1_{\sigma, f}$.
\end{remark}

Later it will be useful to have an expression of $\tilde X_f,
\tilde H_{\sigma, f}$ and $\tilde H^1_{\sigma, f}$ in local
coordinates. Let then $(U,z)$ be a local chart belonging to a
(comfortable, if necessary) atlas adapted to~$S$, so that
$\{\de_1\}$ is a local frame for~$N_S$, and
$\{(\omega^1)^{\otimes\nu_f}\otimes\de_1\}$ is a local frame for
$(N_S^{\otimes\nu_f})^*\otimes N_S$. There is a
holomorphic function~$M_f$ such that
$$
\tilde X_f(\de_1)(\de_{1}^{\otimes\nu_f})=M_f \de_1.
$$
Now, recalling  \pref{eqLfdif}, we obtain
\begin{eqnarray*}
\tilde X_f(\de_1)(\de_{1}^{\otimes\nu_f})&=&\pi\left(\left.\left[\ca
X_f\left(({\de\over
    \de z^1})^{\otimes\nu_f}\right),{\de\over\de z_1}\right]\right|_S\right)\\[8pt]
  &  =&\pi\left(\left.\left[g^j{\de\over\de z^j},{\de\over\de
z^1}\right]\right|_S\right)
    =-\left.{\de g^1\over\de z^1}\right|_S\de_1,
\end{eqnarray*}
and so
\begin{equation}
M_f=-\left.{\de g^1\over\de z^1}\right|_S.
\end{equation}
In particular, recalling that $f$ is tangential we can write $g^1=z^1 h^1$,
and hence
 (5.5) yields
\begin{equation}
M_f=-h^1|_S.
\end{equation}
Similarly, if we write
$\tilde H_{\sigma,f}(\de_1)(\de_{1}^{\otimes\nu_f})= M_{\sigma,f}
\de_1$ and $\tilde H^1_{\sigma,f}(\de_1)(\de_{1})= M^1_{\sigma,f}
\de_1$, we obtain
\begin{equation}
M_{\sigma, f}=M^1_{\sigma, f}=-h^1|_S,
\end{equation}
where $h^1$ is defined by $f^1-z^1=b^1(z^1)^{\nu_f}+h^1 (z^1)^{\nu_f+1}$.

Following ideas originally due to Baum and Bott (see [BB]), we can also
introduce a
holomorphic  action on
the virtual bundle~$TS-N_S^{\otimes \nu_f}$. But let us first define what we
mean by a
holomorphic
action on such a bundle.

\begin{defi}  Let $S^0$ be an open dense subset of a complex manifold $S$, $F$ a
vector bundle
on~$S$,
$X\in H^0(S,\ca T_S\otimes\ca F^*)$, $W$ a vector bundle over~$S^0$ and
$\tilde W$ any
extension
of~$W$ over~$S$ in~$K$-theory. Then we say that $F$ {\it acts
holomorphically on~$\tilde W$}
along~$X$ if $F|_{S^0}$ acts holomorphically on~$W$ along~$X|_{S^0}$.
\end{defi}

Let $S$ be a codimension-one submanifold of~$M$ and take $f\in\End(M,S)$,
$f\not\equiv\id_M$, as
usual. If $f$ is tangential set $X=X_f$. If not, assume that $S$ is
comfortably embedded in
$M$ and
set $X=H_{\sigma,f}$ or $X=H_{\sigma,f}^1$ according to the value
of~$\nu_f$; in this case,
we shall
also assume that $X\not\equiv O$. Set $S^0=S\setminus\Sing(X)$, and let
$\ca Q_f=\ca T_S / X(\ca N_S^{\otimes \nu_f})$. The sheaf $\ca Q_f$
is a coherent analytic sheaf which is locally free over $S^0$. The
associated vector bundle (over $S^0$) is denoted by $Q_f$ and it
is called the {\it normal bundle of $f$}. Then the virtual bundle
$TS-N_S^{\otimes \nu_f}$,
represented by the sheaf~$\ca Q_f$, is an extension (in
the sense of $K$-theory) of~$Q_f$.

\begin{defi}  A {\it Baum-Bott action}  on~$S$ is a holomorphic
action of~$N_S^{\otimes\nu}$ on the virtual bundle
$TS-N_S^{\otimes \nu}$ along a section $X$ of $\ca T_S\otimes
N_S^{\otimes\nu}$, for a suitable $\nu\geq 1$. \pagebreak
\end{defi}

\begin{theorem} \label{dtre} Let $f\in\End(M,S)${\rm ,} $f\not\equiv\id_M${\rm ,} be given.
Assume
that $S$ has codimension one in $M${\rm ,} and that either $f$ is
tangential to $S$ {\rm (}\/and then set $X=X_f${\rm )} or $S$ is comfortably embedded into
$M$ {\rm (}\/and then
set
$X=H_{\sigma,f}$ or $X=H_{\sigma,f}^1$ according to the value of~$\nu_f${\rm ).}
Assume moreover
that
$X \not\equiv 0$. Then there exists a Baum\/{\rm -}\/Bott action
$\tilde B\colon\ca Q_f\to(\ca N_S^{\otimes\nu_f})^*\otimes\ca Q_f$
of $N_S^{\otimes \nu_f}$ on $TS-N_S^{\otimes \nu_f}$ along $X$ defined by
\begin{equation}
\tilde B(s)(u) =\pi_f([X(u),\tilde s])
\end{equation}
where $\pi_f\colon \ca T_S \to \ca Q_f$ is the natural projection{\rm ,} and
$\tilde s\in\ca T_S$
is any
section such that $\pi_f(\tilde s)=s$.
\end{theorem}

{\it Proof.} If $\hat s\in\ca T_S$ is another section such that $\pi_f(\hat s)=s$ we
have $\hat
s-\tilde
s\in X(\ca N_S^{\otimes\nu_f})$; hence $\pi_f([X(u),\hat s-\tilde s])=O$,
and  (5.8) does
not
depend on the choice of~$\tilde s$. Finally, one can easily check that $\tilde B$  is a holomorphic
action on $S^0$.\hfill\qed 

\begin{remark}  Since $S$ has codimension one, $X\colon N_S^{\otimes\nu_f}\to TS$
yields a (possibly
singular) holomorphic foliation on~$S$, and the previous action coincides
with the one used
in~[BB]
for the case of foliations.
\end{remark}

We can also define a third holomorphic action, on the virtual bundle
$TM|_S-N_S^{\otimes\nu_f}$.
Assume that $f$ is tangential, and let
$S^0=S\setminus\Sing(X_f)$, as before. Then the sheaf $\ca W_f=\ca
T_{M,S}/X_f(\ca
N_S^{\otimes\nu_f})$ is a coherent analytic sheaf, locally free over~$S^0$;
let
$W_f=TM|_{S^0}/X_f(N_S^{\otimes\nu_f}|_{S^0})$ be the  associated
vector bundle over~$S^0$. Then the virtual bundle
$TM|_S-N_S^{\otimes \nu_f}$, represented by the sheaf~$\ca W_f$,
is an extension (in the sense of $K$-theory) of~$W_f$.

\begin{defi}  A {\it Lehmann-Suwa action} on~$S$ is a holomorphic
action of~$N_S^{\otimes\nu}$ on~$TM|_S-N_S^{\otimes \nu}$ along a
section~$X$ of $\ca T_S\otimes N_S^{\otimes\nu}$, for a suitable
$\nu\ge 1$.
\end{defi}

Again, the name is chosen to honor the ones who first discovered
the analogous action for holomorphic foliations in any dimension;
see [LS], [LS2] (and [KS] for dimension two).

To present an example of such an action we first need a definition.

\begin{defi}  Let $S$ be a codimension-one, comfortably embedded submanifold of
$M$,
and choose a comfortable atlas~$\gt U$ adapted to a splitting
morphism~$\sigma\colon N_S \to TM|_S$. If $v\in(\ca
N_S^{\otimes\nu})_p$ and $(U,\phe)\in\gt U$ is a chart about~$p\in S$,
we can write $v=\lambda(z'')\de_1^{\otimes\nu}$ for a suitable
$\lambda\in\ca O(U\cap S)$. Then the {\it local extension of~$v$
along the fibers of~$\sigma$} is the local section $\tilde
v=\lambda(z'')(\de/\de z^1)^{\otimes\nu}\in(\ca T_M|_S^{\otimes\nu})_p$.
\end{defi}

If $(\hat U,\hat z)$ is another chart in~$\gt U$ about~$p$, and
$v\in(\ca N_S^{\otimes\nu})_p$, we can also write
$v=\hat\lambda\hat\de_1^{\otimes\nu}$, and we clearly have
$\hat\lambda=(a|_S)^\nu\lambda$. But since $S$ is comfortably
embedded in~$M$ we also have
$$
\left.{\de(\hat\lambda-a^\nu\lambda)\over\de z^1}\right|_S\equiv
0,
$$
and thus
$$
a^\nu\lambda=\hat\lambda+R_2.
$$
Therefore if $\hat v$ denotes the local extension of~$v$ along the
fibers of~$\sigma$ in the chart~$(\hat U,\hat\phe)$ we have
\begin{equation}
\hat v=\hat \lambda\left({\de \over \de
\hz^1}\right)^{\otimes\nu}=a^\nu\lambda {\de
z^{h_1}\over\de\hz^1}\cdots{\de
z^{h_\nu}\over\de\hz^1}\,{\de\over\de
z^{h_1}}\otimes\cdots\otimes{\de\over\de
z^{h_\nu}}+R_2=\tilde{v}+T_1+R_2,
\end{equation}
where
$$
T_1=a\lambda\sum_{\ell=1}^\nu{\de
z^{p_\ell}\over\de\hz^1}{\de\over\de z^1}\otimes\cdots\otimes
{\de\over\de z^{p_\ell}}\otimes\cdots\otimes{\de\over\de z^1}.
$$
Hence:

\begin{theorem} \label{ddue} Let $f\in\End(M,S)${\rm ,} $f\not\equiv\id_M${\rm ,} be given.
Assume
that $S$ is of codimension one and comfortably embedded in $M${\rm ,}
and that $f$ is tangential with $\nu_f>1$. Let $\rho_f\colon \ca T_{M,S} \to
\ca W_f$ be
the
natural projection. Then   a Lehmann\/{\rm -}\/Suwa action $\tilde
V\colon\ca W_f\to(\ca
N_S^{\otimes\nu_f})^*\otimes\ca W_f$ of $N_S^{\otimes \nu_f}$ on
$TM|_S-N_S^{\otimes \nu_f}$ may be defined along $X_f$ by setting 
\begin{equation}
\tilde V(s)(v) =\rho_f([\ca X_f(\tilde v),\tilde s]|_S),
\label{eqV}
\end{equation}
for $s\in\ca W_f$ and $v\in\ca N_S^{\otimes\nu}${\rm ,} where $\tilde s$ is
any element in $\ca T_M|_S$ such that $\rho_f(\tilde s|_S)=s${\rm ,} and
$\tilde v\in\ca T_M|_S^{\otimes\nu_f}$ is an extension of $v$
constant along the fibers of a splitting morphism~$\sigma$.
\end{theorem}

{\it Proof}.  Since $\ca X_f(\tilde v)|_S \in \ca T_{S}$ then
clearly~\pref{eqV} does not depend on the extension $\tilde s$ chosen.
Using~ (5.9) and~(4.7), since $f$ tangential implies $\ca X_f=\ca
H_{\sigma,f}$ and
$T_1^1=R_2$, we have
$$
[\hat{\ca X}_f(\hat v),\tilde s]=[(\ca
X_f+T_1^o+R_2)(\tilde v+T_1+R_2),\tilde s]=[\ca
X_f(\tilde v),\tilde s]+R_1,
$$
and therefore~\pref{eqV} does not depend on the comfortable coordinates
chosen to define it. Finally, arguing as in Theorem~\ref{duno} we can
show that $\tilde V$ actually is a holomorphic action, and we are
done. \hfill\qed

\section{Index theorems for hypersurfaces}

Let $S$ be a compact, globally irreducible, possibly singular hypersurface
in a complex manifold~$M$, and set $S'=S\setminus\Sing(S)$.  Given
the following data:
\begin{itemize}
\ritem{(a)} a line bundle $F$ over $S'$;
\ritem{(b)}a holomorphic section $X$ of $TS'\otimes F^*$;
\ritem{(c)}a vector bundle~$E$ defined on~$M$;
\ritem{(d)}a holomorphic action $\tilde X$ of $F|_{S'}$ on~$E|_{S'}$
along~$X$;
\end{itemize}
\noindent  we can recover a partial connection (in the sense
of Bott) on $E$
restricted to $S^0=S'\setminus\Sing(X)$ as follows: since, by definition
of~$S^0$, the dual
map
$X^*\colon\Xi_X^*\to F^*|_{S^0}$ is an isomorphism, we can define a partial
connection (in
the sense
of Bott~[Bo])
$D\colon\Xi_X\times H^0(S^0,E|_{S^0})\to H^0(S^0,E|_{S^0})$ by setting
$$
D_v(s)=(X^*\otimes\id)^{-1}\bigl(\tilde X(s)\bigr)(v)
$$
for $p\in S^0$, $v\in(\Xi_X)_p$ and $s\in H^0(S^0,E|_{S^0})$. Furthermore,
we can always
extend this partial connection~$D$ to a $(1,0)$-connection on $E|_{S^0}$,
for instance by using a partition of unity (see, e.g., [BB]). Any such connection (which is
a {\it $\Xi_X$-connection} in the terminology of~[Bo], [Su]) will be said to be  {\it induced} by the
holomorphic action~$\tilde X$.

We can then apply the general theory developed by Lehmann and Suwa for
foliations (see in particular Theorem~$1'$ and Proposition~4 of~[LS], as
well as~[Su, Th.~V\negthinspace I.4.8])
to get the following:

\begin{theorem} \label{generale} Let $S$ be a compact{\rm ,} globally irreducible{\rm ,}
possibly singular
hypersurface in an $n$\/{\rm -}\/dimensional complex manifold~$M${\rm ,} and set
$S'=S\setminus\Sing(S)$. Let $F$ be a line bundle over~$S'$
admitting an extension to~$M${\rm ,} and $X$ a holomorphic section of
$TS'\otimes F^*$. Set $S^0=S'\setminus\Sing(X)${\rm ,} and let
$\Sing(S)\cup\Sing(X)=\bigcup_\lambda\Sigma_\lambda$ be the
decomposition of~$\Sing(S)\cup\Sing(X)$ in connected components.
Finally{\rm ,} let $E$ be a vector bundle defined on~$M$. Then for any
holomorphic action~$\tilde X$ of~$F|_{S'}$ on~$E|_{S'}$ along~$X$
and any homogeneous symmetric polynomial~$\phe$ of degree~$n-1${\rm ,} there are  complex numbers $\hbox{\rm
Res}_\phe(\tilde X,E,\Sigma_\lambda)\in\C${\rm ,} depending only on the local behavior
of~$\tilde X$ and~$E$ near~$\Sigma_\lambda${\rm ,} such that
$$
\sum_\lambda\hbox{\rm Res}_\phe(\tilde X,E,\Sigma_\lambda)=\int_S\phe(E),
$$
where $\phe(E)$ is the evaluation of~$\phe$ on the Chern classes of~$E$.
\end{theorem}

Recalling the results of the previous section, we then get the following
index
theorem for holomorphic self-maps:

\begin{theorem} \label{indiceuno} Let $S$ be a compact{\rm ,} globally irreducible{\rm ,}
possibly singular
hypersurface in an $n$\/{\rm -}\/dimensional complex manifold $M$. Let $f
\in \End(M, S)${\rm ,} $f\not\equiv \id_M${\rm ,} be given. Assume that
\begin{itemize}
\ritem{\rm (a)} $f$ is tangential to $S${\rm ,} and   $X=X_f${\rm ,} or that
\ritem{\rm(b)} $S^0=S\setminus\bigl(\Sing(S)\cup
\Sing(f)\bigr)$ is comfortably embedded into $M${\rm ,} and 
$X=H_{\sigma,f}$ if
$\nu_f>1${\rm ,} or
$X=H_{\sigma,f}^1$ if $\nu_f=1$.
\end{itemize}
\noindent Assume moreover $X \not\equiv O$. Let
$\Sing(S)\cup\Sing(X)=\bigcup_\lambda \Sigma_\lambda$ be the decomposition
of
$\Sing(S)\cup\Sing(X)$ in connected components. Finally{\rm ,} let $[S]$ be the
line bundle on~$M$
associated to the divisor~$S$. Then there exist complex numbers~$\hbox{\rm
Res}(X,S,
\Sigma_\lambda)\in\C${\rm ,} depending only on the local
behavior of~$X$ and~$[S]$ near~$\Sigma_\lambda${\rm ,} such that
$$
\sum_{\lambda}\hbox{\rm Res}(X, S,\Sigma_\lambda)=\int_S c_1^{n-1}([S]).
$$
\end{theorem}

\Proof By Theorem~\ref{duno} we have a Camacho-Sad action on~$S$ along $X$ on~$N_{S^0}$. Since
$[S]$ is an
extension to~$M$ of~$N_{S^0}$, we can apply Theorem~\ref{generale}.\hfill\qed 

\begin{remark}  If $M$ has dimension two, and $S$ has
at least one singularity or $X_f$ has at least one zero, then
$S'\setminus \Sing(f)$ is {\it always} comfortably embedded
in $M$. Indeed, it is an open Riemann surface; so $H^1(S'\setminus
\Sing(f), \ca F)=O$ for any coherent analytic sheaf $\ca F$, and the result
follows
from Proposition~\ref{splitting}  and Theorem~\ref{ce}.
\end{remark}

In a similar way, applying [Su, Th.~I\negthinspace V.5.6],
Theorem~\ref{ddue}, and recalling
that
$\phe(H-L)=\phe(H\otimes L^*)$ for any vector bundle $H$, line bundle~$L$
and homogeneous
symmetric
polynomial~$\phe$, we get

\begin{theorem} \label{lesuvariation} Let $S$ be a compact{\rm ,} globally irreducible{\rm ,}
possibly singular hypersurface in an $n$\/{\rm -}\/dimensional complex
manifold $M$. Let $f \in\End(M,S)${\rm ,} $f \not\equiv \id_M${\rm ,} be
given. Assume that $S'=S\setminus\Sing(S)$ is comfortably embedded
into $M${\rm ,} and that $f$ is tangential to $S$ with $\nu_f>1$. Let
$\Sing(S)\cup\Sing(X_f)=\bigcup_\lambda \Sigma_\lambda$ be the
decomposition of $\Sing(S)\cup\Sing(X_f)$ in connected components.
Finally{\rm ,} let $[S]$ be the line bundle on $M$ associated to the
divisor $S$. Then for any homogeneous symmetric polynomial~$\phe$
of degree~$n-1$ there exist complex numbers $\hbox{\rm
Res}_\phe(X_f,TM|_S-[S]^{\otimes\nu_f},\Sigma_\lambda)\in \C${\rm ,}
depending only on the local behavior of~$X_f$
and~$TM|_S-[S]^{\otimes \nu_f}$ near~$\Sigma_\lambda${\rm ,} such that
$$
\sum_\lambda\hbox{\rm Res}_\phe(X_f,TM|_S-[S]^{\otimes
\nu_f},\Sigma_\lambda)= \int_S
\varphi\bigl(TM|_S\otimes([S]^*)^{\otimes \nu_f}\bigr).
$$
\end{theorem}

Finally, applying the Baum-Bott index theorem (see  [Su,
Th.~I\negthinspace
I\negthinspace I.7.6]) and Theorem~5.2 we get

\begin{theorem} \label{baumbot} Let $S$ be a compact{\rm ,} globally irreducible{\rm ,}
 smooth complex
hypersurface
in an $n$\/{\rm -}\/dimensional complex manifold $M$. Let $f\in \End(M,S)${\rm ,}
$f \not\equiv \id_M${\rm ,} be given. Assume that
\begin{itemize}
\ritem{\rm (a)} $f$ is tangential to $S${\rm ,} and  $X=X_f${\rm ,} or that
\ritem{\rm(b)} $S^0=S\setminus \Sing(f)$  is comfortably embedded into $M${\rm ,} and 
$X=H_{\sigma,f}$ if
$\nu_f>1${\rm ,} or
$X=H_{\sigma,f}^1$ if $\nu_f=1$.
\end{itemize}
\noindent Assume moreover $X \not\equiv O$. Let
$\Sing(X)=\bigcup_\lambda \Sigma_\lambda$ be the
decomposition of $\Sing(X)$ in connected components.
Finally{\rm ,} let $[S]$ be the line bundle on $M$ associated to the
divisor $S$. Then for any homogeneous symmetric polynomial~$\phe$
of degree~$n-1$ there exist complex numbers $\hbox{\rm
Res}_\phe(X,TS-[S]^{\otimes \nu_f},\Sigma_\lambda)\in \C${\rm ,}
depending only on the local behavior of~$X$ and~$TS-[S]^{\otimes
\nu_f}$ near~$\Sigma_\lambda${\rm ,} such that
$$
\sum_\lambda\hbox{\rm Res}_\phe(X,TS-[S]^{\otimes
\nu_f},\Sigma_\lambda)=\int_S \phe\bigl(TS\otimes([S]^*)^{\otimes
\nu_f}\bigr).
$$
\end{theorem}
\vglue12pt

Thus, we have recovered three main index theorems of foliation theory in the
setting of
holomorphic self-maps fixing pointwise a hypersurface.

Clearly, these index theorems are as useful as   the
formulas for the computation of the residues are explicit; the rest of this
section is devoted to deriving such formulas in many important
cases.

Let us first describe the general way these residues are defined
in Lehmann-Suwa theory. Assume   the hypotheses of
Theorem~\ref{generale}. Let $\tilde U_0$ be a tubular neighborhood
of~$S^0$ in~$M$, and denote by $\rho\colon\tilde U_0\to S^0$ the
associated retraction. Given any connection~$D$ on~$E|_{S^0}$
induced by the holomorphic action~$\tilde X$ of~$F$ along~$X$, set
$D^0=\rho^*(D)$. Next, choose an open set $\tilde U_\lambda\subset
M$ such that $\tilde
U_\lambda\cap\bigl(\Sing(S)\cup\Sing(X)\bigr)=\Sigma_\lambda$, and
a compact real $2n$-dimensional manifold $\tilde
R_\lambda\subset\tilde U_\lambda$ with $C^\infty$ boundary
containing~$\Sigma_\lambda$ in its interior and such that
$\de\tilde R_\lambda$ intersects $S$ transversally. Let
$D^\lambda$ be any connection on~$E|_{\tilde U_\lambda}$, and
denote by $B\bigl(\phe(D^0),\phe(D^\lambda)\bigr)$ the Bott
difference form of~$\phe(D^0)$ and~$\phe(D^\lambda)$ on $\tilde
U_0\cap\tilde U_\lambda$. Then (see [LS] and [Su, Chap.\
I\negthinspace V])
\begin{equation}
\hbox{\rm Res}_\phe(\tilde
X,E,\Sigma_\lambda)=\int_{R_\lambda}\phe(D^\lambda)
-\int_{\de R_\lambda} B\bigl(\phe(D^0),\phe(D^\lambda)\bigr),
\label{residuouno}
\end{equation}
where $R_\lambda=\tilde R_\lambda\cap S$. A similar formula holds for
virtual
vector bundles too; see again~[Su, Chap.~I\negthinspace V].

\begin{remark}  When $\Sigma_\lambda=\{x_\lambda\}$ is an isolated singularity
of~$S$, the second
integral in
\pref{residuouno} is taken on the link of~$x_\lambda$ in $S$. In particular if
$S$ is
not irreducible at $x_k$ then the residue is the sum of several
terms, one for each irreducible component of $S$ at $x_k$.
\end{remark}

We now specialize \pref{residuouno} to our situation. Let us begin
with the Camacho-Sad action: we shall compute the residues for
connected components~$\Sigma_\lambda$ reduced to an isolated point
$x_\lambda$.  Let again $[S]$ be the line bundle associated to the
divisor~$S$, and choose an open set $\tilde U_\lambda\subset M$
containing $x_\lambda$ so that $\tilde
U_\lambda\cap\bigl(\Sing(S)\cup \Sing(X)\bigr)=\{x_\lambda\}$ and
$[S]$ is trivial on $\tilde U_\lambda$; take  as $D^\lambda$
the trivial connection for~$[S]$ on $W$ with respect to some
frame. In particular, then, $\phe(D^\lambda)=O$ on $R_\lambda$. By
\pref{residuouno} the residue \pagebreak is then obtained simply by integrating
$B\bigl(\phe(D^0),\phe(D^\lambda)\bigr)$ over $\de R_\lambda$.
Notice furthermore that since $[S]$ is a line bundle there is only
one nontrivial~$\phe$ to consider: the $(n-1)^{\rm th}$ power of the
linear symmetric function, so that~$\phe(D)=c_1^{n-1}([S])$.

Let $\eta^j$ be a connection one-form of $D^j$, for
$j=0$,~$\lambda$; with respect to a suitable frame for~$[S]$ we
can assume that $\eta^\lambda\equiv O$. Let
$$
\tilde{\eta}:=t\eta^0+(1-t)\eta^\lambda=t \eta^0,
$$
and let $\tilde{K}:=d \tilde \eta+ \tilde \eta \wedge \tilde \eta
= d \tilde \eta$. From the very definition of the Bott difference
form, it follows that
$$
B\bigl(\phe(D^0),\phe(D^\lambda)\bigr)= \left({1\over2\pi
i}\right)^{n-1}\int_0^1
\tilde{K}^{n-1}.
$$
A straightforward computation shows that
$$
\tilde{K}^{n-1}=(n-1)t^{n-2}dt \wedge \eta^0 \wedge
\overbrace{d\eta^0 \wedge \cdots \wedge d \eta^0}^{n-2}+N,
$$
where $N$ is a term not containing $dt$. Therefore
\begin{equation}
B\bigl(\phe(D^0),\phe(D^\lambda)\bigr)=\left({1\over2\pi
i}\right)^{n-1}\eta^0 \wedge
\overbrace{d\eta^0 \wedge \cdots \wedge d
\eta^0}^{n-2}.\label{bottdiff}
\end{equation}

Assume now that $x_\lambda \in \Sing(X)$ and $S$ is smooth at
$x_\lambda$. Up to shrinking $\tilde U_\lambda$ we may assume that
$\tilde U_\lambda$ is the domain of a chart~$z$ adapted to $S$
(and belonging to a comfortable atlas if necessary), so that
$\{\de_1\}$ is a local frame for $N_{S^0}$, and
$\{dz^2,\ldots,dz^n\}$ is a local frame for~$T^*S^0$. Then any
connection $D$ induced by the Camacho-Sad action is locally
represented by the (1,0)-form $\eta^0$ such that
$D(\de_1)=\eta^0\otimes\de_1$. To compute~$\eta^0$, we first of
all notice that $X=g^p{\de\over\de z^p}\otimes (\omega^1)^{\otimes
{\nu_f}}$, if $X=X_f$ or $X=H_{\sigma, f}$, and
$X=(1+b^1)g^p{\de\over\de z^p}\otimes \omega^1$ if~$X=H^1_{\sigma,
f}$. Then, when $X$ is $X_f$ or $H_{\sigma, f}$,
$$
(X^*)^{-1}((\omega^1)^{\otimes\nu_f})={1\over g^p}\,dz^p,
$$
while when $X=H^1_{\sigma, f}$,
$$
(X^*)^{-1}((\omega^1)^{\otimes\nu_f})={1\over(1+b^1)g^p}\,dz^p.
$$
Therefore, recalling formulas  (5.6) and  (5.7), we can choose $D$ so
that when $X$ is
$X_f$ or
$H_{\sigma, f}$,
\begin{equation}
\eta^0=(X^*\otimes\id)^{-1}\bigl(\tilde X(\de_1)\bigr)=-\left.{h^1\over
g^p}\right|_S \,
dz^p,\label{connectionuno}
\end{equation}
while when $X=H^1_{\sigma, f}$,
\begin{equation}
\eta^0=(X^*\otimes\id)^{-1}\bigl(\tilde H^1_{\sigma, f}
(\de_1)\bigr)=-\left.{h^1\over (1+b^1)g^p}\right|_S \,
dz^p.\label{connectionunouno}
\end{equation}

\begin{remark}  When $n=2$ and $X=X_f$ we recover the connection form obtained
in [Br]. The form $\eta$ introduced in~[A2], which is the opposite
of~$\eta^0$, is the connection form of the dual connection
on~$N_{S^0}^*$, by [A2, (1.7)]. Since the definition of Chern
class implicitly used in~[A2] is the opposite of the one used
in~[Br] everything is coherent. Finally, when $n=2$ and
$X=H^1_{\sigma,f}$ we have obtained the correct multiple of the
form~$\eta$ introduced in~[A2] when $S$ was the smooth zero
section of a line bundle (notice that $1+b^1$ is constant because
$S$ is compact, and that the form $\eta$ of [A2] must be divided
by~$b=1+b^1$ to get a connection form).
\end{remark}

Now we can take $R_1=\{|g^p(x)|\le\eps \mid p=2,\ldots,n\}$ for a
suitable $\eps>0$ small enough.
In particular, if we set
$\Gamma=\{|g^p(x)|=\eps\mid p=2,\ldots,n\}\cap S$, oriented so
that $d\theta^2\wedge \cdots \wedge d\theta^n
>0$ where $\theta^p=\arg(g^p)$, then arguing as in\break [L, \S 5] or [LS,
\S 4] (see
also
[Su,~pp.105--107]) we obtain
\begin{equation}
\hbox{\rm Res}(X,S, \{x_\lambda\})=\left({-i\over 2\pi
}\right)^{n-1}\int_{\Gamma}{(h^1)^{n-1} \over g^2 \cdots
g^n}\,dz^2 \wedge \cdots \wedge dz^n,
\label{residuoespluno}
\end{equation}
when $X=X_f$ or $H_{\sigma,f}$, while when $X=H^1_{\sigma, f}$ we
have
\begin{equation}
\hbox{\rm Res}(H^1_{\sigma, f}, S,\{x_\lambda\})=\left({-i\over 2\pi
}\right)^{n-1}\int_{\Gamma}{(h^1)^{n-1} \over (1+b^1)^{n-1}g^2
\cdots g^n}\,dz^2 \wedge \cdots \wedge dz^n.
\label{residuoespldue}
\end{equation}

\begin{remark}  For $n=2$,  formulas \pref{residuoespluno} and \pref{residuoespldue} give the
indices defined in~[A2]. Thus, if $S$ is smooth,
Theorem~\ref{indiceuno} implies the index theorem of~[A2], because
$c_1([S])=c_1(N_S)$. In an analogous way, Lehmann and Suwa (see
[L], [LS], [LS2]) proved that the Camacho-Sad index theorem also is a
consequence of Theorem~\ref{generale}.
\end{remark}

When $x_\lambda$ is an isolated singular point of $S$ the
computation of the residue is more complicated, because one cannot apply
directly the results in [LS] as before, for in general there is no
natural extension of $\Xi_X$ and the Camacho-Sad action to
$\Sing(S)$. However we are able to compute explicitly the index
in this case too when~$n=2$, and when $n>2$ and $f$ is tangential with
$\nu_f>1$.

If $n=2$ we can choose local coordinates $\{(w^1, w^2)\}$ in $\tilde
U_\lambda$
so that $S \cap \tilde U_\lambda=\{l(w^1, w^2)=0\}$ for some holomorphic
function
$l$, and $dl \wedge dw^2 \neq 0$ on $S \cap \tilde U_\lambda\setminus \{
x_\lambda\}$. In particular $(l,w^2)$ are local
coordinates adapted to~$S^0$ near $S\cap\tilde U_\lambda\setminus \{
x_\lambda\}$ and
${\de \over \de l}$ can be chosen as a local frame for
$N_{S^0}$ on $\de R_1$.

\begin{remark}  When $S^0$ is comfortably embedded in~$M$ the
chart $(l,w^2)$ should belong to a comfortable atlas. Studying the
proofs of Propositions~\ref{splitting}  and of Theorem~\ref{ce}\ one sees
that this is possible up to replacing~$l$ by a function of the
form $\hat l=\bigl(1+c(w^2)l)l$, where $c$ is a holomorphic
function defined on~$S\cap\tilde U_\lambda\setminus \{
x_\lambda\}$. Since to \pagebreak compute the residues we only need the
behavior of $l$ and $w^2$ near~$\de R_1$, it is easy to check that
using $\hat l$ or~$l$ in the following computations yields the same
results. So for the sake of simplicity we shall not distinguish
between $l$ and~$\hat l$ in the sequel.
\end{remark}

Up to shrinking $\tilde U_\lambda$, we can again assume that $[S]$
is trivial on~$\tilde U_\lambda$. The function $l$ is a local
generator of $\ca I_S$ on $\tilde U_\lambda$. Then the dual of
$[l]\in \ca I_S/\ca I_S^2$, denoted by $s$, is a holomorphic frame
of $[S]$ on $\tilde U_\lambda$ which extends the holomorphic frame
$\de \over \de l$ of $N_{S'}$ (see [Su,~p.86]). In particular
 $s|_{\de R_1}={\de \over \de l}$. We then
choose on~$[S]|_{\tilde U_\lambda}$ the trivial connection with
respect to~$s$, so that~$\eta^\lambda=O$. We are left with the
computation of the form $\eta^0$ near $\de R^1$. But if $X=X_f$
or $X=H_{\sigma,f}$ we can apply~\pref{connectionuno} to get
$$
\eta^0|_{\de R_1}=-\left.{(l \circ f - l)-b^1l^{\nu_f}\over l\cdot(w^2 \circ
f -
w^2)}\right|_{\de
R_1}dw^2,
$$
where
$$
b^1=\left.{l\circ f-l\over l^{\nu_f}}\right|_S
$$
is identically zero when $f$ is tangential. On the other hand, when
$X=H^1_{\sigma,f}$,
applying
\pref{connectionunouno}  we get
$$
\eta^0|_{\de R_1}=-\left.{(l \circ f - l)-b^1l\over (l+(l\circ f - l))(w^2
\circ f - w^2)}
\right|_{\de R_1}dw^2.
$$
Hence  the residue is
\begin{equation}
\hbox{\rm Res}(X, S, \{x_\lambda\})={1\over 2 \pi i}\int_{\de
R_1}\left.{(l \circ f - l)-b^1l^{\nu_f}\over l\cdot(w^2 \circ f -
w^2)}\right|_{S}dw^2,
\label{residuodue}
\end{equation}
when $X=X_f$ or $X=H_{\sigma,f}$, while when $X=H^1_{\sigma,f}$,
\begin{equation}
\hbox{\rm Res}(H^1_{\sigma,f}, S, \{x_\lambda\})={1\over 2 \pi
i}\int_{\partial R_1}\left.{(l \circ f - l)-b^1l\over (l+(l\circ f - l))(w^2
\circ f -
w^2)}\right|_{S}dw^2.
\label{residuodueuno}
\end{equation}
\vglue12pt

\begin{remark}  When $f$ is tangential we have $b^1\equiv0$; therefore the formula
\pref{residuodue}
gives the
index defined in~[BT], and Theorem~\ref{indiceuno} implies the index
theorem of~[BT].
\end{remark}

When $n>2$, $f$ is tangential and $\nu_f>1$, we can define a local vector
field $\tilde v_f$ which generates the Camacho-Sad action $\tilde
X_f$ and compute explicitly the residue even at a singular point
$x_\lambda$ of $S$. To see this, assume $(w^1,\ldots, w^n)$ are
local coordinates in $\tilde U_\lambda$ so that $S \cap \tilde
U_\lambda=\{l(w^1,\ldots,
w^n)=0\}$ for some holomorphic function $l$. Define the vector
field $\tilde{v}_f$ on $\tilde U_\lambda$ by
\begin{equation}
\tilde{v}_f={w^1\circ f - w^1\over l^{\nu_f}}{\de\over\de
w^1}+\ldots+{w^n\circ f - w^n\over l^{\nu_f}}{\de\over\de
w^n}.
\label{vecfield}
\end{equation}
We claim that the ``holomorphic action'' $\theta_{\tilde{v}_f}$ in
the sense of Bott [Bo]  of $\tilde{v}_f$ on $N_{S'}$ as defined
in~[LS,~p.177] coincides with our Camacho-Sad action, and thus we
can apply [LS, Th.~1] to compute the residue. To prove this we
consider $W_1=\{x \in \tilde U_\lambda| {\de l \over \de w^1}(x) \neq 0\}$.
On
this open set we make the following change of coordinates:
$$
\left\{ \begin{array}{ll}
z^1=l(w^1,\ldots, w^n),\\  z^p=w^p& \hbox{for $p=2,\ldots,n$}. \end{array}\right.
$$
The new coordinates $(z^1,\ldots, z^n)$ are adapted to $S$ on~$W_1$.
If $f^j=z^j+g^j(z^1)^{\nu_f}$ as
usual, we have
\begin{equation}
w^p \circ f - w^p=g^p (z^1)^{\nu_f},
\label{primop}
\end{equation}
and
\begin{equation}
w^1 \circ f - w^1= {\de w^1\over \de z^j}g^j
(z^1)^{\nu_f}+R_{2\nu_f}=\left({\de l\over\de
w^1}\right)^{-1}\left[g^1-{\de l\over \de
w^p}g^p\right](z^1)^{\nu_f}+R_{2\nu_f}.
\label{primouno}
\end{equation}
Therefore, from \pref{primop} and \pref{primouno}, taking into account that
 $\nu_f>1$, we get
\begin{eqnarray}
\tilde{v}_f&=&\left({w^1 \circ f - w^1\over (z^1)^{\nu_f}}{\de
l\over \de w^1}+{w^p \circ f - w^p\over (z^1)^{\nu_f}}{\de l\over
\de w^p}\right){\de\over \de z^1}\\
&&+{w^q \circ f - w^q\over
(z^1)^{\nu_f}}{\de\over\de z^q}=\ca X_f(\de_1^{\otimes
\nu_f})+R_2,\label{vtilde}\nonumber
\end{eqnarray}
which gives the claim on $W_1$. Since the same holds on each
$W_j=\break \{x \in \tilde U_\lambda | {\de l\over \de w^j}(x)\neq 0\}$,
$j=1,\ldots,n$,
and $(\tilde U_\lambda\cap S)\setminus\{x_\lambda\}=\bigcup_j W_j$, it
follows that the Bott
holomorphic action induced by $\tilde{v}_f$ is the same as the Camacho-Sad
action
given by $\tilde X_f$. Thus, if we choose --- as we can --- the
coordinates $(w^1,\ldots, w^n)$ as in~[LS, Th.~2], that is so that
$\{l, (w^p \circ f
- w^p)/
l^{\nu_f}\}$ form a regular sequence at $x_\lambda$, the residue is
expressed by
the formula after~[LS,~Th.~2]. Taking into account that, since $f$ is tangential and
by
\pref{vtilde}, the function $l$ divides $dl(\tilde{v}_f)$, we get
\begin{equation}
\hbox{\rm Res}(X_f, S, \{x_\lambda\})=\left({-i\over 2\pi
i}\right)^{n-1} \int_\Gamma {\left[\sum_{j=1}^n{\de l\over \de
w^j}(w^j \circ f - w^j) \right]^{n-1}\over l^{n-1}\prod_{p=2}^n
(w^p \circ f - w^p)}\,dw^2 \wedge \cdots \wedge dw^n,
\label{resnugrande}
\end{equation}
where this time
$$
\Gamma=\left\{w\in\tilde U_\lambda \biggm| \left|{w^p \circ f - w^p\over
l^{\nu_f}}(w)\right|=\epsilon,\ l(w)=0\right\},
$$
for a suitable $0<\epsilon<<1$, and $\Gamma$ is oriented
as usual (in particular $\Gamma=(-1)^{[{n-1\over 2}]} R_{u_0}$
where $R_{u_0}$ is the set defined in [LS,~Th.~2]).

Note that for $n=2$ we recover, when $\nu_f>1$,
formula~\pref{residuodue}. On the other hand, if $x_\lambda$ is nonsingular
for $S$, then the previous argument with $l=w^1$ works  for
$\nu_f=1$ as well, and we get formula~\pref{residuoespluno}.

Summing up, we have proved the following:

\begin{theorem} \label{calcoloresuno} Let $S$ be a compact{\rm ,} globally irreducible{\rm ,}
possibly
singular
hypersurface in an $n$\/{\rm -}\/dimensional complex manifold $M$. Let $f
\in \End(M, S)${\rm ,} $f\not\equiv \id_M${\rm ,} be given. Assume that
\begin{itemize}
\ritem{\rm (a)} $f$ is tangential to $S${\rm ,} and $X=X_f${\rm ,} or that
\ritem{\rm(b)} $S^0=S\setminus\bigl(\Sing(S)\cup
\Sing(f)\bigr)$ is comfortably embedded into $M${\rm ,} and 
$X=H_{\sigma,f}$ if
$\nu_f>1${\rm ,} or
$X=H_{\sigma,f}^1$ if $\nu_f=1$.
\end{itemize}
\noindent Assume $X \not\equiv O$. Let $x_\lambda\in S$ be an isolated point
of~$\Sing(S)\cup\Sing(X)$. Then the number~$\hbox{\rm Res}(X,S,
\{x_\lambda\})\in\C$ introduced in Theorem~{\rm \ref{indiceuno}} is given
\begin{itemize}
\ritem{\rm(i)} if $x_\lambda\in\Sing(X)\cap (S\setminus \Sing (S))${\rm ,} and $f$ is tangential 
or $S^0$ is
comfortably
embedded
in~$M$ and $\nu_f>1${\rm ,} by
$$
\hbox{\rm Res}(X,S, \{x_\lambda\})=\left({-i\over 2\pi
}\right)^{n-1}\int_{\Gamma}{(h^1)^{n-1} \over g^2 \cdots
g^n}\,dz^2 \wedge \cdots \wedge dz^n;
$$
\ritem{\rm(ii)} if $x_\lambda\in\Sing(X)\cap (S\setminus \Sing (S))${\rm ,} $S^0$ is comfortably 
embedded
in~$M$ and $\nu_f=1${\rm ,} by
$$
\hbox{\rm Res}(H^1_{\sigma, f}, S,\{x_\lambda\})=\left({-i\over 2\pi
}\right)^{n-1}\int_{\Gamma}{(h^1)^{n-1} \over (1+b^1)^{n-1}g^2
\cdots g^n}\,dz^2 \wedge \cdots \wedge dz^n;
$$
\ritem{\rm(iii)} if $n=2${\rm ,} $x_\lambda\in\Sing(S)${\rm ,} and $f$ is tangential or
$S^0$ is
comfortably
embedded in~$M$ and $\nu_f>1${\rm ,} by
$$
\hbox{\rm Res}(X, S, \{x_\lambda\})={1\over 2 \pi i}\int_{\de
R_1}\left.{(l \circ f - l)-b^1l^{\nu_f}\over l\cdot(w^2 \circ f -
w^2)}\right|_{S}dw^2;
$$
\ritem{\rm(iv)} if $n=2${\rm ,} $x_\lambda\in\Sing(S)${\rm ,} $S^0$ is comfortably
embedded
in~$M$ and $\nu_f=1${\rm ,} by
$$
\hbox{\rm Res}(H^1_{\sigma,f}, S, \{x_\lambda\})={1\over 2 \pi
i}\int_{\partial R_1}\left.{(l \circ f - l)-b^1l\over (l+(l\circ f - l))(w^2
\circ f -
w^2)}\right|_{S}dw^2;
$$
\ritem{\rm(v)} if $n>2${\rm ,} $x_\lambda\in\Sing(S)${\rm ,} $f$ is tangential and
$\nu_f>1${\rm ,} by
\end{itemize}
$$
\hbox{\rm Res}(X_f, S, \{x_\lambda\})=\left({-i\over 2\pi
i}\right)^{n-1} \int_\Gamma {\left[\sum_{j=1}^n{\de l\over \de
w^j}(w^j \circ f - w^j) \right]^{n-1}\over l^{n-1}\prod_{p=2}^n
(w^p \circ f - w^p)}\,dw^2 \wedge \cdots \wedge dw^n.
$$
\end{theorem}
\vglue12pt

Our next aim is to compute the residue for the Lehmann-Suwa
action, at least for an isolated smooth
point~$x_\lambda\in\Sing(X_f)$. Let $(W,w)$ be a local chart
about~$x_\lambda$ belonging to a comfortable atlas. Set $l=w^1$
and define $\tilde{v}_f$ as in \pref{vecfield}. By \pref{vtilde} the
Lehmann-Suwa action $\tilde V$ is given by the holomorphic action
(in the sense of Bott) of~$\tilde{v}_f$ on $TM|_S-[S]^{\otimes
\nu_f}$. Therefore we can apply [L], [LS] (see also [Su,~Ths.\
I\negthinspace V.5.3, I\negthinspace V.5.6], and [Su,
Remark~I\negthinspace V.5.7]) to obtain
$$
\hbox{\rm Res}_\phe(X_f,TM|_S-[S]^{\otimes \nu_f},
\{x_\lambda\})=\hbox{Res}_\phe(X_f,TM|_S,\{x_\lambda\}),
$$
where $\hbox{\rm Res}_\phe(X_f,TM|_S,\{x_\lambda\})$ is the
residue for the local Lie derivative action of $\tilde{v}_f$ on
$TM|_S$ given by
$$
\tilde V_l(s)(\tilde{v}_f)= [\tilde{v}_f, \tilde s]|_S,
$$
where $s$ is a section of $TM|_S$ and $\tilde s$ is a local
extension of $s$ constant along the fibers of $\sigma$.

We can write an expression of $\tilde V_l$ in local coordinates. Let
$(U,z)$ be a local chart belonging to a comfortable atlas. Then
$\{{\de\over\de
z^1},\ldots, {\de \over \de z^n}\}$ is a local frame for $TM$, and
$\{(\omega^1)^{\otimes \nu_f}\otimes\left.{\de\over\de z^1}\right|_S,\ldots,
(\omega^1)^{\otimes \nu_f}\otimes\left. {\de \over \de z^n}\right|_S\}$ is a
local frame for $(N_S^{\otimes\nu_f})^\ast\otimes TM|_S$. Thus
there exist holomorphic functions $V_j^k$ (for
$j$,~$k=1,\ldots,n$) so that
$$
\tilde V_l({\de\over\de z^j})(\de_1^{\otimes \nu_f})=V_j^k{\de
\over \de z^k}.
$$
Now, from \pref{hestende} we get
\begin{eqnarray*}
\tilde V_l({\de\over\de z^j})(\de_1^{\otimes
\nu_f})&=&\left.\left[\ca X_f\left(({\de\over \de
z^1})^{\otimes \nu_f} \right), {\de\over\de
z^j}\right]\right|_S\\
&=&\left.\left[h^1z^1{\de\over\de z^1}+g^p{\de
\over\de z^p}, {\de \over\de z^j}
\right]\right|_S=-h^1|_S\delta_j^1{\de \over \de z^1} - \left.{\de
g^p \over \de z^j}\right|_S {\de \over \de z^p},
\end{eqnarray*}
and hence
\begin{equation}
V_1^1= -h^1|_S,\qquad  V_p^1\equiv 0,\qquad  V_j^p=-\left.{\de g^p
\over \de z^j}\right|_S.
\label{connev}
\end{equation}
Therefore [Su, Th.~I\negthinspace V.5.3] yields

\begin{theorem} \label{calcoloresdue} Let $S$ be a compact{\rm ,} globally irreducible{\rm ,}
possibly singular hypersurface in an $n$\/{\rm -}\/dimensional complex
manifold $M$. Let $f \in\End(M,S)${\rm ,} $f \not\equiv \id_M${\rm ,} be
given. Assume that $S'=S\setminus\Sing(S)$ is comfortably embedded
into~$M${\rm ,} and that $f$ is tangential to $S$ with $\nu_f>1$. Let
$x_\lambda\in\Sing(X_f)$ be an isolated smooth point of
$\Sing(S)\cup\Sing(X_f)$. Then for any homogeneous symmetric
polynomial~$\phe$ of degree~$n-1$ the complex number $$\hbox{\rm
Res}_\phe(X_f,TM|_S-[S]^{\otimes\nu_f},\{x_\lambda\})$$ introduced
by Theorem~{\rm \ref{lesuvariation}} is given by
\begin{equation}
\hbox{\rm Res}_\phe(X_f,TM|_S-[S]^{\otimes \nu_f},
\{x_\lambda\})=\int_\Gamma{\varphi(V)\,dz^2\wedge \cdots \wedge
dz^n\over g^2\cdots g^n},
\label{eqrls}
\end{equation}
where $V=(V_j^k)$ is the matrix given by {\rm \pref{connev}} and $\Gamma$ is
as in {\rm (6.5)}.
\end{theorem}

\begin{remark}  We adopt here the convention that if $V$ is an $n\times n$ matrix
then $c_j(V)$ is the $j^{\rm th}$ symmetric function of the eigenvalues $V$
multiplied
by $(i/2\pi)^{j}$, and $\phe(V)=\phe\bigl(c_1(V),\ldots,c_{n-1}(V)\bigr)$.
\end{remark}

Finally, if $x_\lambda$ is an isolated point
in~$\Sing(X)$, the complex numbers\break $\hbox{\rm
Res}_\phe(X,TS-[S]^{\otimes \nu_f},\{x_\lambda\})$ appearing in
Theorem~\ref{baumbot} can be computed exactly as in the foliation case
using a Grothendieck residue with a formula very similar
to~\pref{eqrls}; see~[BB], [Su, Th.~I\negthinspace I\negthinspace
I.5.5].

\section{Index theorems in higher codimension}

Let $S\subset M$ be a complex submanifold of codimension $1<m<n$ in a
complex
$n$-manifold~$M$. A
way to get index theorems for holomorphic self-maps of~$M$ fixing
pointwise~$S$ is to
blow-up $S$
and then apply the index theorems for hypersurfaces; this is what we plan to
do in this
section.

We shall denote by $\pi\colon M_S\to M$ the blow-up of~$M$
along~$S$, and by $E_S=\pi^{-1}(S)$ the exceptional divisor, which
is a  hypersurface in~$M_S$ isomorphic to the projectivized normal
bundle~$\P(N_S)$.

\begin{remark}  If $S$ is singular, the blow-up $M_S$ is in general singular too. So
this approach
works
only for smooth submanifolds.
\end{remark}

If $(U,z)$ is a chart adapted to~$S$ centered in~$p\in S$, in $M_S$ we have
$m$ charts $(\tilde
U_r,w_r)$
centered in $[\de_1],\ldots,[\de_m]$ respectively, where if $v\in N_{S,p}$,
$v\ne
O$, we denote by~$[v]$ its projection in~$\P(N_S)$. The coordinates $z^j$
and $w^h_r$ are
related by
$$
z^j(w_r)=\left\{\begin{array}{ll}w^j_r& \hbox{if $j=r,m+1,\ldots,n$,}\\ 
w^r_r w^j_r& \hbox{if $j=1,\ldots,r-1,r+1,\ldots,m$.}\end{array}\right.
$$

\begin{remark}  We have $\tilde U_r\cap E_S=\{w^r_r=0\}$, and thus $(\tilde
U_r,w_r)$ is adapted
to~$E_S$ up
to a permutation of the coordinates.
\end{remark}

Now take $f\in\End(M,S)$, $f\not\equiv\id_M$, and assume that $df$ acts
as the identity on~$N_S$ (this is automatic if~$\nu_f>1$, while if
$\nu_f=1$ it happens if and only if $f$ is tangential).  Then we
can lift~$f$ to a unique map $\tilde f\in\End(M_S,E_S)$, $\tilde
f\not\equiv\id_{M_S}$, such that $f\circ\pi=\pi\circ\tilde f$ (see, e.g.,
[A1] for details). If $(U,z)$ is a chart adapted to~$S$ and we set
$f^j=z^j\circ f$ and $\tilde f^j_r=w^j_r\circ\tilde f$, 
\begin{equation}
\tilde f^j_r(w_r)=\left\{\begin{array}{ll} {\displaystyle f^j\bigl(z(w_r)\bigr)}& \hbox{if
$j=r,m+1,\ldots,n$,}\\ 
\noalign{\smallskip}
{\displaystyle {f^j\bigl(z(w_r)\bigr)\over f^r\bigl(z(w_r)\bigr)}}&\hbox{if
$j=1,\ldots,r-1,r+1,\ldots,m$.}\end{array}\right.
\label{eqstuno}
\end{equation}
If $f$ is tangential we can find holomorphic functions
$h^r_{r_1\ldots r_{\nu_f+1}}$ symmetric in the lower indices such
that
\begin{equation}
f^r-z^r=h^r_{r_1\ldots r_{\nu_f+1}}
    z^{r_1}\cdots z^{r_{\nu_f+1}}+R_{\nu_f+2};
\label{eqstunomez}
\end{equation}
as usual, only the restriction to~$S$ of each $h^r_{r_1\ldots r_{\nu_f+1}}$
is uniquely
defined. Set
then
$$
Y=h^r_{r_1\ldots
r_{\nu_f+1}}|_S\,\de_r\otimes\omega^{r_1}\otimes\cdots\otimes\omega^{r_{\nu_
f+1}};
$$
it is a local section of $N_S\otimes(N_S^*)^{\otimes(\nu_f+1)}$.

On the other hand, if $f$ is not tangential set
$B=(\pi\otimes\id)_*\circ X_f$, where $\pi\colon TM|_S\to N_S$ is
the canonical projection. In this way we get a global section of
$N_S\otimes(N_S^*)^{\otimes\nu_f}$, not identically zero if and
only if $f$ is not tangential, and given in local adapted
coordinates by
$$
B=g^r_{r_1\ldots
r_{\nu_f}}|_S\,\de_r\otimes\omega^{r_1}\otimes\cdots\otimes\omega^{r_{\nu_f}
}.
$$

\begin{defi}  Take $p\in S$. If $f$ is tangential, a non-zero vector
$v\in(N_S)_p$ is a {\it singular direction} for~$f$ at~$p$ if
$X_f(v\otimes\cdots\otimes
v)=O$ and
$Y(v\otimes\cdots\otimes v)\wedge v=O$. If $f$ is not tangential, $v$ is a
{\it singular
direction}
if $B(v\otimes\cdots\otimes v)\wedge v=O$.
\end{defi}

\begin{remark}  The condition $Y(v\otimes\cdots\otimes v)\wedge v=O$ is equivalent
to requiring
$Y(v\otimes\cdots\otimes v)=\lambda v$ for some $\lambda\in\C$.
\end{remark}

Of course, in the tangential case we must check that this definition is
well-posed, because
the
morphism~$Y$ depends on the local coordinates chosen to define it. First of
all, if $(U,z)$
is a
chart adapted to~$S$ and centered at~$p$ then $X_f(v\otimes\cdots\otimes
v)=O$ when $f$ is
tangential means
\begin{equation}
g^p_{r_1\ldots r_{\nu_f}}(O)\,v^{r_1}\cdots v^{r_{\nu_f}}\,{\de\over\de
z^p}=O,
\label{eqstdue}
\end{equation}
where $v=v^r\de_r$. Now let $(\hat U,\hat z)$ be another chart adapted
to~$S$ centered
in~$p$. Then
we can find holomorphic functions
$a^r_s$ and $\hat a^r_s$ such that $\hat z^r=a^r_s z^s$ and~$z^r=\hat
a^r_s\hat z^s$.
Arguing as in
the proof of~\pref{eqtanggen} we get
$$
a^{r_1}_{s_1}\cdots a^{r_{\nu_f+1}}_{s_{\nu_f+1}}\hat h^r_{r_1\ldots
r_{\nu_f+1}}=
a^r_s h^s_{s_1\ldots s_{\nu_f+1}}+\sum_{\ell=1}^{\nu_f+1}{\de
a^r_{s_\ell}\over\de z^p}
g^p_{s_1\ldots\hat{s_\ell}\ldots s_{\nu_f+1}}+R_1,
$$
where the index with the hat is missing from the list. Therefore
$$
\hat Y=Y+\hat a^s_r\sum_{\ell=1}^{\nu_f+1}\left.{\de a^r_{s_\ell}\over\de
z^p}
g^p_{s_1\ldots\hat{s_\ell}\ldots
s_{\nu_f+1}}\right|_S\,\de_s\otimes\omega^{s_1}\otimes\cdots
\otimes\omega^{s_{\nu_f+1}};
$$
in particular if $X_f(v\otimes\cdots\otimes v)=O$ equation \pref{eqstdue} yields
$$
\hat Y(v\otimes\cdots\otimes v)=Y(v\otimes\cdots\otimes v),
$$
and the notion of singular direction when $f$ is tangential is well-defined.

\begin{prop} \label{stuno} Let $S\subset M$ be a complex submanifold of
codimension
$1<m<n$ of a
complex $n$\/{\rm -}\/manifold~$M${\rm ,} and take $f\in\End(M,S)${\rm ,} $f\not\equiv\id_M${\rm ,} such that
$df$ acts as the
identity
on~$N_S$ {\rm (}\/that is $f$ is tangential{\rm ,} or $\nu_f>1${\rm ,} or both\/{\rm ).} Denote by
$\pi\colon
M_S\to M$ the blow\/{\rm -}\/up of~$M$ along~$S$ with exceptional divisor~$E_S${\rm ,} and
let $\tilde
f\in\End(M_S,E_S)$ be the lifted map. Then\/{\rm :}
\begin{itemize}
\ritem{\rm (i)}if $S$ is comfortably embedded in $M$ then $E_S$ is
comfortably embedded
in~$M_S${\rm ,} and
the choice of a splitting morphism for~$S$ induces a splitting morphism
for~$E_S${\rm ;}
\ritem{\rm(ii)} $d\tilde f$ acts as the identity on~$N_{E_S}${\rm ;}
\ritem{\rm(iii)} $\tilde f$ is always tangential{\rm ;} furthermore $\nu_{\tilde
f}=\nu_f$ if $f$
is
tangential{\rm ,} $\nu_{\tilde f}=\nu_f-1$ otherwise\/{\rm ;}
\ritem{\rm(iv)} a direction $[v]\in E_S$ is a singular point for~$\tilde f$
if and only if it is a
singular  direction
for~$f$.
\end{itemize}
\end{prop}

\Proof (i) Let $\gt U=\{(U_\alpha,z_\alpha)\}$ be a comfortable atlas adapted
to~$S$; we claim
that
$\tilde{\gt U}=\{(\tilde U_{\alpha,r},w_{\alpha,r})\}$ is a comfortable
atlas adapted
to~$E_S$ (and
in particular determines a splitting morphism for~$E_S$). Let us first prove
that it is a
splitting
atlas, that is that
$$
\left.{\de w^j_{\beta,s}\over\de w^r_{\alpha,r}}\right|_{E_S}\equiv 0
$$
for every $r$, $s$, $j\ne s$ and indices $\alpha$ and $\beta$. We
have
$$
z^j_\beta=z^j_\beta|_S+\left.{\de z^j_\beta\over\de z^s_\alpha}\right|_S
z^s_\alpha+
{1\over 2}\left.{\de^2 z^j_\beta\over\de z^u_\alpha\de z^v_\alpha}\right|_S
z^u_\alpha
z^v_\alpha+R_3.
$$
Since $w^r_{\alpha,r}=z^r_\alpha$, if $j=p>m$ we immediately get
$$
\left.{\de w^p_{\beta,s}\over\de w^r_{\alpha,r}}\right|_{E_S}=\left.{\de
z^p_\beta\over\de
z^r_\alpha}\right|_S\equiv 0,
$$
because $\gt U$ is a splitting atlas. If $j=t\le m$,
\begin{eqnarray}
z^t_\beta&=&\left.{\de z^t_\beta\over\de z^s_\alpha}\right|_S z^s_\alpha+
{1\over 2}\left.{\de^2 z^t_\beta\over\de z^u_\alpha\de z^v_\alpha}\right|_S
z^u_\alpha
z^v_\alpha+R_3\label{eqsttre}
\\
&=&\left[\left.{\de z^t_\beta\over\de z^r_\alpha}\right|_S
+\sum_{u\ne r}
\left.{\de z^t_\beta\over\de z^u_\alpha}\right|_S
w^u_{\alpha,r}\right]w^r_{\alpha, r}+
O\bigl((w^r_{\alpha,r})^3\bigr),
\nonumber
\end{eqnarray}
because $\gt U$ is a comfortable atlas. Therefore if $t\ne s$,
$$
w^t_{\beta,s}={z^t_\beta\over z^s_\beta}={\left.{\de z^t_\beta\over\de
z^r_\alpha}\right|_S
+\sum_{u\ne r}\left.{\de z^t_\beta\over\de z^u_\alpha}\right|_S
w^u_{\alpha,r}+O\bigl((w^r_{\alpha,r})^2\bigr)\over
\left.{\de z^s_\beta\over\de z^r_\alpha}\right|_S
+\sum_{u\ne r}\left.{\de z^s_\beta\over\de z^u_\alpha}\right|_S
w^u_{\alpha,r}+O\bigl((w^r_{\alpha,r})^2\bigr)},
$$
and so
$$
{\de w^t_{\beta,s}\over\de w^r_{\alpha,r}}=O(w^r_{\alpha,r}),
$$
as required.

Finally, since $w^s_{\beta,s}=z^s_\beta$, equation \pref{eqsttre} yields
$$
{\de^2 w^s_{\beta,s}\over\de(w^r_{\alpha,r})^2}=O(w^r_{\alpha,r}),
$$
and $\tilde{\gt U}$ is a comfortable atlas, as claimed.
\vglue1pt
(ii) Since $df$ acts as the identity on $N_S$, in local coordinates we can
write
$$
f^j(z)=z^j+g^j_{r_1\ldots r_{\nu_f+1}}z^{r_1}\cdots
z^{r_{\nu_f}}+R_{\nu_f+1},
$$
with $g^s_{r_1}|_S\equiv 0$ if $\nu_f=1$. Then \pref{eqstuno} yields
\begin{equation}
\tilde f^j_r(w_r)=w^j_r+(w^r_r)^{\nu_f} g^j_{r_1\ldots
r_{\nu_f}}\bigl(z(w_r)\bigr)
    w_r^{\hat r_1}\cdots w_r^{\hat r_{\nu_f}}+O\bigl((w_r^r)^{\nu_f+1})
\label{eqstquattro}
\end{equation}
if $j=r,m+1,\ldots,n$, and
\begin{eqnarray}\qquad 
\tilde f^j_r(w_r)&=&w^j_r+(w_r^r)^{\nu_f-1}\bigl[g^j_{r_1\ldots
r_{\nu_f}}\bigl(z(w_r)\bigr)-w^j_r
g^r_{r_1\ldots r_{\nu_f}}\bigl(z(w_r)\bigr)\bigr]w_r^{\hat r_1}\cdots
w_r^{\hat
r_{\nu_f}}\label{eqstcinque}\\
&&+O\bigl((w_r^r)^{\nu_f})
\nonumber
\end{eqnarray}
if $j=1,\ldots,r-1,r+1,\ldots,m$, where $w^{\hat s}_r=w^s_r$ if $s\ne r$,
and $w^{\hat
r}_r=1$. In
particular, $d\tilde f$ acts as the identity on~$N_{E_S}$.

\vskip1pt
(iii) We have
$$
g^j_{r_1\ldots r_{\nu_f}}|_{E_S}\bigl(z(w_r)\bigr)=g^j_{r_1\ldots
r_{\nu_f}}|_S(O,w''_r);
$$
therefore if $f$ is tangential then $w^r_r$ divides all $g^s_{r_1\ldots
r_{\nu_f}}\bigl(z(w_r)\bigr)$,
while it does not divide some $g^p_{r_1\ldots r_{\nu_f}}\bigl(z(w_r)\bigr)$.
In particular,
then,
$\tilde f$ is tangential and $\nu_{\tilde f}=\nu_f$, by \pref{eqstquattro} and
\pref{eqstcinque}. On
the other
hand, if $f$ is not tangential $w^r_r$ does not divide some $g^s_{r_1\ldots
r_{\nu_f}}\bigl(z(w_r)\bigr)$; therefore
\begin{eqnarray*}
&&\hskip-36pt
\bigl[g^s_{r_1\ldots r_{\nu_f}}\bigl(z(w_r)\bigr)-w^s_r
g^r_{r_1\ldots r_{\nu_f}}\bigl(z(w_r)\bigr)\bigr]\bigr|_{E_S}\\
&&\qquad\quad =g^s_{r_1\ldots
r_{\nu_f}}(O,w''_r)
-w^s_rg^r_{r_1\ldots r_{\nu_f}}(O,w''_r)\not\equiv 0,
\end{eqnarray*}
and thus $\nu_{\tilde f}=\nu_f-1$ and $\tilde f$ is again tangential.
\vskip1pt

(iv) Take $v\in(N_S)_p$, $v\ne O$, and a chart $(U,z)$ adapted
to~$S$ centered in~$p$. Then $v=v^s\de_s$, with $v^r\ne 0$ for
some $r$. Hence $[v]\in\tilde U_r$ has coordinates
$$
w^j_r([v])=\left\{\begin{array}{ll}
0&\hbox{if $j=r,m+1,\ldots,n$,}\\ 
v^j/v^r& \hbox{if $j=1,\ldots,r-1,r+1,\ldots,m$.}
\end{array}\right.
$$
If $f$ is not tangential, then $[v]$ is a singular point
for~$\tilde f$ if and only if
$$
[v^rg^s_{r_1\ldots r_{\nu_f}}(O)-v^sg^r_{r_1\ldots
r_{\nu_f}}(O)]v^{r_1}\cdots v^{r_{\nu_f}}=0
$$
for all $s$, \pagebreak and thus if and only if $B(v\otimes\cdots\otimes
v)\wedge v=O$, as claimed. 

If $f$ is tangential, writing $f^s-z^s$ as in \pref{eqstunomez} we get
\begin{eqnarray*} 
\tilde f^s_r(w_r)&=&w^s_r+(w_r^r)^{\nu_f}\bigl[h^s_{r_1\ldots
r_{\nu_f+1}}\bigl(z(w_r)\bigr)-w^s_r
h^r_{r_1\ldots r_{\nu_f+1}}\bigl(z(w_r)\bigr)\bigr]w_r^{\hat r_1}\cdots
w_r^{\hat
r_{\nu_f+1}}\\
&&+O\bigl((w_r^r)^{\nu_f+1})
\end{eqnarray*}
for $s\ne r$, and then it is clear that $[v]$ is a singular point
for~$\tilde f$ if and only if $v$ is a singular direction
for~$f$.\Endproof 

We therefore get index theorems in any codimension:

\begin{theorem} \label{stdue} Let $S$ be a compact complex submanifold of
codimension $1<m<n$ in
an $n$\/{\rm -}\/dimensional complex manifold $M$. Let $f
\in \End(M,S)${\rm ,} $f\not\equiv \id_M${\rm ,} be given, and assume that $df$ acts as
the identity
on~$N_S$.
Let $\bigcup_\lambda \Sigma_\lambda$ be the decomposition in connected
components of the
set of
singular directions for~$f$ in~$\P(N_S)$. Then there exist complex
numbers~$\hbox{\rm
Res}(f,S,
\Sigma_\lambda)\in\C${\rm ,} depending only on the local
behavior of~$f$ and~$S$ near~$\Sigma_\lambda${\rm ,} such that
$$
\sum_{\lambda}\hbox{\rm Res}(f,
S,\Sigma_\lambda)=\int_{E_S}c_1^{n-1}([E_S])=\int_S\pi_*
c_1^{n-1}([E_S]),
$$
where $E_S$ is the exceptional divisor in the blow-up~$\pi\colon M_S\to M$
of $M$ along
$S${\rm ,} and
$\pi_*$ denotes the integration along the fibers of the
bundle~$\pi|_{E_S}\colon E_S\to S$.
\end{theorem}

\Proof  This follows immediately from Theorem~\ref{indiceuno}, Proposition~\ref{stuno}, and
the projection
formula (see, e.g., [Su,~Prop.~I\negthinspace I.4.5]).\hfill\qed 

\begin{remark}  The restriction to $E_S$ of the cohomology class $c_1([E_S])$ is the
Chern class
$\zeta=c_1(T)$
of the tautological bundle $T$ on the bundle $\pi|_{E_S}\colon E_S\to S$
and
it
satisfies the relation
\begin{eqnarray*}
 \zeta^{n-m}-\pi|_{E_S}^*c_1(N_S)\zeta^{n-m-1}&+&\pi|_{E_S}^*c_2(N_S)\zeta^{n-m
-2}+
\cdots\\
&&\cdots +(-1)^{n-m}\pi|_{E_S}^*c_{n-m}(N_S)=0
\end{eqnarray*}
in the cohomology ring of the bundle (see, e.g.,  [GH, pp.~606--608]). This
formula can
sometimes be used to compute $\zeta$ in terms of the Chern classes of~$N_S$
and~$TM$ in
specific examples.
\end{remark}

\begin{theorem} \label{sttre} Let $S$ be a compact complex submanifold of
codimension $1<m<n$ in
an $n$\/{\rm -}\/dimensional complex manifold $M$. Let $f
\in\End(M,S)${\rm ,} $f \not\equiv \id_M${\rm ,} be given{\rm ,} and set $\nu=\nu_f$ if $f$ is
tangential{\rm ,} and
$\nu=\nu_f-1$ otherwise. Assume that $S$ is comfortably embedded into $M${\rm ,}
and that
$\nu>1$. Let
$\bigcup_\lambda \Sigma_\lambda$ be the decomposition in connected
components of the set of
singular directions for~$f$ in~$\P(N_S)$. Finally{\rm ,} let $\pi\colon M_S\to M$
be the blow\/{\rm -}\/up
of~$M$ along~$S${\rm ,} with exceptional divisor~$E_S$. Then for any homogeneous
symmetric
polynomial~$\phe$ of degree~$n-1$ there exist complex numbers $\hbox{\rm
Res}_\phe(f,TM_S|_{E_S}-N_{E_S}^{\otimes\nu},\Sigma_\lambda)\in
\C${\rm ,} depending only on the local behavior of~$f$
and~$TM_S|_{E_S}-N_{E_S}^{\otimes\nu}$
near~$\Sigma_\lambda${\rm ,} such that
$$
\sum_\lambda\hbox{\rm
Res}_\phe(f,TM_S|_{E_S}-N_{E_S}^{\otimes\nu},\Sigma_\lambda)=
\int_S\pi_* \varphi\bigl(TM_S|_{E_S}\otimes(N_{E_S}^*)^{\otimes\nu}\bigr),
$$
where $\pi_*$ denotes the integration along the fibers of the bundle $E_S\to
S$.
\end{theorem}

\Proof This  follows immediately from Theorem~\ref{lesuvariation}, Proposition~\ref{stuno},
and the
projection
formula.\hfill\qed 

\begin{theorem} \label{stquattro} Let $S$ be a compact complex submanifold of
codimension $1<m<n$
in
an $n$\/{\rm -}\/dimensional complex manifold $M$. Let $f\in \End(M,S)${\rm ,}
$f \not\equiv \id_M${\rm ,} be given{\rm ,} and assume that $df$ acts as the identity
on~$N_S$. Set
$\nu=\nu_f$
if $f$ is tangential{\rm ,} and $\nu=\nu_f-1$ otherwise. Let
$\bigcup_\lambda \Sigma_\lambda$ be the decomposition in connected
components of the set of
singular directions for~$f$ in~$\P(N_S)$. Finally{\rm ,} let $\pi\colon M_S\to M$
be the blow\/{\rm -}\/up
of~$M$ along~$S${\rm ,} with exceptional divisor~$E_S$. Then for any homogeneous
symmetric
polynomial~$\phe$ of degree~$n-1$ there exist complex numbers
$\hbox{\rm Res}_\phe(f,TE_S-N_{E_S}^{\otimes \nu},\Sigma_\lambda)\in \C${\rm ,}
depending only on
the
local behavior of~$f$ and~$TE_S-N_{E_S}^{\otimes \nu}$
near~$\Sigma_\lambda${\rm ,} such
that
$$
\sum_\lambda\hbox{\rm Res}_\phe(f,TE_S-N_{E_S}^{\otimes
\nu},\Sigma_\lambda)=\int_S\pi_*
\phe\bigl(TE_S\otimes(N_{E_S}^*)^{\otimes \nu}\bigr),
$$
where $\pi_*$ denotes the integration along the fibers of the bundle $E_S\to
S$.
\end{theorem}

\Proof This  follows immediately from Theorem~\ref{baumbot}, Proposition~\ref{stuno}, and
the projection
formula.\hfill\qed

\section{Applications to dynamics}

We conclude this paper with   applications to the study of the
dynamics of
endomorphisms of complex manifolds,  first recalling  a definition from [A2]:

\begin{defi}  Let $f\in\End(M,p)$ be a germ at~$p\in M$ of a holomorphic self-map of
a complex
manifold~$M$ fixing $p$. A {\it parabolic curve} for~$f$ at~$p$ is
a injective holomorphic map $\phe\colon\Delta\to M$ satisfying the
following properties:
\begin{itemize}
\ritem{(i)} $\Delta$ is a simply connected domain in~$\C$ with
$0\in\de\Delta$; 
\ritem{(ii)} $\phe$ is continuous at the origin,
and $\phe(0)=p$; \ritem{(iii)} $\phe(\Delta)$ is invariant
under~$f$, and $(f|_{\phe(\Delta)})^n\to p$ as $n\to\infty$.
\end{itemize}
\noindent Furthermore, we say that $\phe$ is {\it tangent to a
direction~$v\in T_pM$} at~$p$ if for one (and hence any) chart
$(U,z)$ centered at~$p$ the direction of
$z\bigl(\phe(\zeta)\bigr)$ converges to the direction $dz_p(v)$
as~$\zeta\to0$.
\end{defi}

Now we have the promised dynamical interpretation of~$X_f$ at nonsingular
points:

\begin{prop} \label{intdyn} Assume that $S$ has codimension one in~$M${\rm ,} and
take $f\in\End(M,S)${\rm ,} $f\not\equiv\id_M$. Let $p\in S$ be a
regular point of~$X_f${\rm ,} that is\break $X_f(p)\ne O$. Then
\begin{itemize}
\ritem{\rm(i)}If $f$ is tangential then no infinite orbit of $f$ can stay
arbitrarily close to~$p$. More precisely{\rm ,} there is a neighborhood
$U$ of~$p$ such that for every $q\in U$ there is $n_0\in\N$ such
that $f^{n_0}(q)\notin U$ or $f^{n_0}(q)\in S$.
\ritem{\rm(ii)}If $\Xi_f(p)$ is transversal to~$T_pS$ {\rm (}\/so in particular $f$
is non-tangential\/{\rm )} and $\nu_f>1$ then there exists at least one
parabolic curve for $f$ at~$p$ tangent to~$\Xi_f(p)$.
\ritem{(iii)}If $\Xi_f(p)$ is transversal to~$T_pS${\rm ,} $\nu_f=1${\rm ,} and
$|b(p)|\ne0${\rm ,}~$1$ or $b(p)=\exp(2\pi i\theta)$ where $\theta$
satisfies the Bryuno  condition {\rm (}\/and $b$ is the
function defined in Remark~{\rm 1.1)} then there exists an $f$\/{\rm -}\/invariant
one\/{\rm -}\/dimensional   holomorphic disk~$\Delta$ passing
through~$p$ tangent to~$\Xi_f(p)$ such that $f|_\Delta$ is
holomorphically conjugated to  the multiplication
by~$b(p)$.
\end{itemize}
\end{prop}

\Proof In local adapted coordinates centered at~$p\in S$ we can write
$$
f^j(z)=z^j+(z^1)^{\nu_f}g^j(z),
$$
so that
$$
\Xi_f(p)=\hbox{\rm Span}\left(g^1(O)\left.{\de\over\de
z^1}\right|_p+\cdots
    +g^n(O)\left.{\de\over\de z^n}\right|_p\right).
$$
In case (i), we have $g^1=z^1h^1$ for a suitable holomorphic
function $h^1$, and $g^{p_0}(O)\ne0$ for some $2\le p_0\le n$,
because $p$ is not singular. Therefore we can apply [AT,
Prop.~2.1] (see also [A2, Prop.~2.1]), and the
assertion follows.

Now, $\Xi_f(p)$ is transversal to~$T_pS$ if and only if $g^1(O)\ne
0$. In case (ii) we can then write
$$
f^j(z)=z^j+g^j(O)(z^1)^{\nu_f}+O(\|z\|^{\nu_f+1})
$$
with $g^1(O)\ne0$. Then $\Xi_f(p)$ is a non-degenerate
characteristic direction of~$f$ at~$p$ in the sense of Hakim, and
thus by [H1,~2] there exist at least $\nu_f-1$ parabolic curves
for $f$ at~$p$ tangent to~$\Xi_f(p)$.

If $\nu_f=1$, it is easy to see that $b^1(p)=1+g^1(O)$, and
$b^1(p)\ne 1$ because $\Xi_f(p)$ is transversal to~$T_pS$.
Therefore we can write
$$
f^j(z)=\left\{\begin{array}{ll} b^1(p)z^1+O(\|z\|^2)&\hbox{if $j=1$,}\\ 
z^j+g^j(O)z^1+O(\|z\|^2)&\hbox{if $2\le j\le n$,}\end{array}\right.
$$
and the assertion in case (iii) follows immediately from
[P\"o] (see also [N]).\Endproof 

In other words, $X_f$ essentially dictates the dynamical behavior
of $f$ away from the singularities --- or, from another point of
view, we can say that the interesting dynamics is concentrated
near the singularities of~$X_f$.

\begin{remark}  If $\Xi_f(p)$ is transversal to $T_pS$, $\nu_f=1$ and $b(p)=0$ or
$b(p)=\exp(2\pi i\theta)$ with $\theta$ irrational not satisfying the
Bryuno condition, there might still be an $f$-invariant
one-dimensional holomorphic disk passing through~$p$ and tangent
to~$\Xi_f(p)$. On the other hand, if $b(p)=\exp(2\pi i\theta)$ is
a $k^{\rm th}$ root of unity, necessarily different from one, several
things might happen. For instance, if $b(p)=-1$, up to a linear
change of coordinates we can write
$$
f^j(z)=\left\{\begin{array}{ll} z^1+z^1\bigl(-2+(z^1)^{\mu_1}\hat g^1(z)\bigr)&\hbox{if
$j=1$,}\\  z^j+(z^1)^{\mu_j+1}\hat g^j(z)&\hbox{if $j=2,\ldots,n$,}\end{array}\right.
$$
for suitable $\mu_1,\ldots,\mu_n\in\N$ and holomorphic functions
$\hat g^j$ not divisible by~$z^1$ and such that $\hat g^j(O)=0$ if
$\mu_j=0$. Then if $\mu_1=0$,
\begin{eqnarray*}
&&(f\circ f)^j(z)\\
&&\qquad  =\left\{\begin{array}{ll} z^1-z^1\bigl[\hat g^1(z)+\hat
g^1\bigl(f(z)\bigr)-\hat g^1(z)
    \hat g\bigl(f(z)\bigr)\bigr]&\hbox{if $j=1$,}\\ 
z^j+(z^1)^{\mu_j+1}\bigl[\hat g^j(z)-\bigl(-1+\hat
g^1(z)\bigr)^{\mu_j+1}\hat g^j\bigl(f(z)\bigr)\bigr]&\hbox{if
$j=2,\ldots,n$.}\end{array}\right.
\end{eqnarray*}
So $\nu_{f\circ f}=1$, $f\circ f$ is non-tangential but $p$ is
singular for $f\circ f$. On the other hand, if $\mu_1=1$,
\begin{eqnarray*}
&&(f\circ f)^j(z)\\
&&\qquad =\left\{\begin{array}{ll} z^1-(z^1)^2\bigl[\hat g^1(z)-\hat g^1\bigl(
    f(z)\bigr)+O(z^1)\bigr]&\hbox{if $j=1$,}\\ 
z^j+(z^1)^{\mu_j+1}\bigl[\hat g^j(z)+(-1)^{\mu_j}\hat
g^j\bigl(f(z)\bigr)+O(z^1)\bigr]&\hbox{if $j=2,\ldots,n$.}\end{array}\right.
\end{eqnarray*}
Now if, for instance, $\mu_2=0$ we get $\nu_{f\circ f}=1$, but
$f\circ f$ is tangential and $p$ is singular for $f\circ f$. But
if $\mu_2=2$ and $\mu_j\ge 2$ for $j\ge 3$ we get $\nu_{f\circ
f}=3$ and $p$ can be either singular or nonsingular for $f\circ
f$.
\end{remark}

\begin{remark}  If $\nu_f=1$, $\Xi_f(p)\subset T_pS$ and $S$ is compact, necessarily
$f$ is
tangential, because $b\equiv 1$ and then $g^1(0,z'')\equiv 0$. If
$S$ is not compact we might have an isolated point of tangency,
and in that case we might have parabolic curves at~$p$ not tangent
to~$\Xi_f(p)$. For instance, the methods of [A1] show that this
happens for the map
$$
f^j(z)=\left\{\begin{array}{ll} z^1+z^1\bigl(az^2+bz^3+h_1(z'')+z^1h_2(z)\bigr)&\hbox{if
$j=1$,}\\  z^2+z^1\bigl(c+h_3(z)\bigr)&\hbox{if $j=2$,}\\ 
z^3+z^1g^3(z)&\hbox{if $j=3$,}\end{array}\right.
$$
when $a$,~$c\ne 0$.
\end{remark}

Finally, we describe a couple of applications to endomorphisms of complex
surfaces:

\begin{Cor} \label{appluno} Let $S$ be a smooth compact one\/{\rm -}\/dimensional
submanifold of
a complex surface~$M${\rm ,} and take $f\in\End(M,S)${\rm ,} $f\not\equiv\id_M$. Assume
that $f$ is
tangential{\rm ,} or that $S\setminus\Sing(f)$ is comfortably embedded in~$M${\rm ,} and
let $X$
denote $X_f${\rm ,} $H_{\sigma,f}$ or $H^1_{\sigma,f}$ as usual\/{\rm ;} assume moreover
that
$X\not\equiv O$. Then
\begin{itemize}
 \ritem{(i)} if $c_1(N_S)\ne 0$ then $\chi(S)-\nu_fc_1(N_S)>0${\rm ;}
\ritem{(ii)} if $c_1(N_S)>0$ then $S$ is rational{\rm ,} $\nu_f=1$ and $c_1(N_S)=1$.
\end{itemize}
\end{Cor}

\Proof The well-known theorem about the localization of the top Chern class at
the zeroes of a
global
section (see, e.g., [Su,~Th.~I\negthinspace I\negthinspace I.3.5])
yields
\begin{equation}
\sum_{x\in\Sing(X)} N(X;x)=\chi(S)-\nu_f\, c_1(N_S),
\label{eqcasosemplice}
\end{equation}
where $N(X;x)$ is the multiplicity of $x$ as a zero of $X$. Now,
If $c_1(N_S)\ne 0$ then by Theorem~\ref{indiceuno} the set $\Sing(X)$ is not
empty.
Therefore the sum in \pref{eqcasosemplice}
must be strictly positive, and the assertions follow.\hfill\qed 

\begin{defi}  Let $f\in\End(M,S)$, $f\not\equiv\id_M$. We say that a point
$p\in S$ is {\it weakly attractive} if there are infinite orbits
arbitrarily close to~$p$, that is, if for every neighborhood $U$
of~$p$ there is $q\in U$ such that $f^n(q)\in U\setminus S$ for
all~$n\in\N$. In particular, this happens if there is an infinite
orbit converging to~$p$.
\end{defi}

Then we can prove the following

\begin{prop} \label{dina} Let $S$ be a smooth compact one\/{\rm -}\/dimensional
submanifold of
a complex surface~$M${\rm ,} and take $f\in\End(M,S)${\rm ,} $f\not\equiv\id_M$. If $f$
is tangential
then there are at most $\chi(S)-\nu_f c_1(N_S)$ weakly attractive
points for $f$ on $S$.
\end{prop}

\Proof By \pref{eqcasosemplice} the sum of zeros of the section $X_f$
(counting multiplicity) is equal to $\chi(S)-\nu_f c_1(N_S)$. Thus
the number of zeros (not counting multiplicity) is at most
$\chi(S)-\nu_f c_1(N_S)$. The assertion then follows from
Proposition~\ref{intdyn}. \phantom{abitcold}\Endproof 

Finally, the previous index theorems allow a
classification of the smooth curves which are fixed by a
holomorphic map and are dynamically trivial.

\begin{theorem} \label{duna} Let $S$ be a smooth compact one\/{\rm -}\/dimensional
submanifold of
a complex surface~$M${\rm ,} and take $f\in\End(M,S)${\rm ,} $f\not\equiv\id_M$.
Moreover
assume that $\sp(df_p)=\{1\}$ for some $p \in S$. If there
are no weakly attractive points for $f$ on $S$ then only one of the
following cases
occurs\/{\rm :}
\begin{itemize}
\ritem{(i)} $\chi(S)=2${\rm ,} $c_1(N_S)=0${\rm ,} or
\ritem{(ii)} $\chi(S)=2${\rm ,} $c_1(N_S)=1${\rm ,} $\nu_f=1${\rm ,} or
\ritem{(iii)} $\chi(S)=0${\rm ,} $c_1(N_S)=0$.
\end{itemize}
\end{theorem}

\Proof Since $N_S$ is a line bundle over a compact curve $S$, the
action of $df$ on $N_S$ is given by multiplication by a constant,
and since $df_p$ has only the eigenvalue $1$ then this constant
must be 1. If $f$ were nontangential then by
Proposition~\ref{intdyn}.(ii) all but a finite number of points of $S$
would be weakly attractive. Therefore $f$ is tangential. By [A2,
Cor.~3.1] (or [Br, Prop.~7.7]) if there is a point $q\in S$
so that $\hbox{Res}(X_f, N_S, p)\not\in \Q^{+}$ then $q$ is weakly
attractive. Thus the sum of the residues
is nonnegative and by Theorem~\ref{indiceuno} it follows that~$c_1(N_S)\geq 0$.
Thus
\pref{eqcasosemplice} yields
\begin{equation}
\chi(S)\geq \nu_f c_1(N_S) \geq 0.
\label{stima}
\end{equation}
Therefore the only possible cases are $\chi(S)=0,2$. If
$\chi(S)=0$ then \pref{stima} implies $c_1(N_S)=0$. Assume that
$\chi(S)=2$. Thus $c_1(N_S)=0$, 1, 2. However if $c_1(N_S)=1$ and
$\nu_f=2$ or if $c_1(N_S)=2$ (and necessarily~$\nu_f=1$) then
\pref{eqcasosemplice} would imply that $X_f$ has no zeroes, and thus
$c_1(N_S)=0$ by Theorem~\ref{indiceuno}.\phantom{whatawaste} \Endproof

\end{document}